\documentclass[a4paper, twoside]{article}
\usepackage[inner=2cm, outer=2cm, top=2.6cm, bottom=2.3cm, headsep=0.7cm]{geometry}
\usepackage{amsmath}
\numberwithin{equation}{section}
\usepackage{amssymb}
\usepackage{mathrsfs}
\usepackage{upgreek}
\usepackage{dsfont}
\usepackage{mathdots}
\usepackage{mathtools}
\usepackage[all]{xy}
\usepackage{tikz}
\usetikzlibrary{calc, graphs, positioning, quotes}
\usepackage{longtable}
\usepackage{color}
\usepackage{graphicx}

\usepackage{fancyhdr}
\usepackage{enumitem}
\setlist[enumerate]{itemsep=3pt, parsep=0pt, topsep=0pt,  label=(\arabic*)}
\setlist[itemize]{itemsep=0pt, parsep=0pt, topsep=0pt}

\usepackage[thmmarks, amsmath]{ntheorem}
\newtheorem{Theorem}{Theorem}
\newtheorem{Lemma}[Theorem]{Lemma}
\newtheorem{Proposition}[Theorem]{Proposition}
\newtheorem{Corollary}[Theorem]{Corollary}

{\theorembodyfont{\normalfont}
\newtheorem{Remark}[Theorem]{Remark}

\newtheorem{Example}[Theorem]{Example}
\newtheorem{Examples}[Theorem]{Examples}
}

{\theoremstyle{empty}
\theoremsymbol{\mbox{$\Box$}}
\theoremheaderfont{\itshape}
\theorembodyfont{\normalfont}
\newtheorem{proof}{}}
\newenvironment{Proof}[1][Proof:]{\begin{proof}[#1]} {\end{proof}}

\numberwithin{Theorem}{section}
\numberwithin{equation}{section}

\usepackage{titlesec}
\titleformat{\section} {\centering\large\bfseries}{\S \thesection}{10pt}{}

\usepackage[numbers, square]{natbib}
\usepackage{hyperref}
\hypersetup{bookmarksopen=true, bookmarksnumbered=true, pdfborder={0 0 0}}

\newcommand\AC{{\mathbb{C}[q^{\pm \frac{1}{2}}]}}
\newcommand\nT{\widetilde{T}}

\newcommand\lrleq{\underset {\mathrm{LR}} {\leqslant}}

\newcommand\lrequiv{\underset{\mathrm{LR}}{\sim}}

\newcommand\hq[1][1]{q^{ \frac{#1}{2}}}
\newcommand\mhq[1][1]{q^{- \frac{#1}{2}}}
\newcommand\af{\boldsymbol{a}}
\newcommand\Vg{V_{\text{geom}}}
\newcommand\dleq{\underset{\mathrm{d}}{\leqslant}}
\newcommand\dequiv{\underset{\mathrm{d}}{\sim}}

\begin{document}
\fancyhead[LE, RO]{\thepage}
\fancyhead[CO]{\scshape Representations of Coxeter groups and homology of Coxeter graphs}
\fancyhead[CE]{\scshape H. Hu}

\title{Representations of Coxeter groups and homology of Coxeter graphs}
\author{Hongsheng Hu}
\date{}
\maketitle

\let\thefootnote\relax\footnotetext{{\bfseries Keywords:} Coxeter groups, IR-representations, homology of graphs}

\let\thefootnote\relax\footnotetext{{\bfseries MSC 2020:} 20C15, 20F55}

\pdfbookmark[1]{Abstract}{abstract}
\abstract{We classify a class of complex representations of an arbitrary Coxeter group via characters of the integral homology of certain graphs. Such representations can be viewed as a generalization of the geometric representation and correspond to the second-highest 2-sided cell in the sense of Kazhdan-Lusztig.
We also give a description of the cell representation provided by this 2-sided cell, and find out all its simple quotients for simply laced Coxeter system with no more than one circuit in the Coxeter graph.}

\setcounter{tocdepth}{1}
\tableofcontents

\large

\section{Introduction}
Let $W = \langle S \mid s^2 = (st)^{m_{st}} = e, \forall s \neq t \in S \rangle$ be a Coxeter group of finite rank.
Weyl groups and affine Weyl groups are important examples, playing crucial roles in Lie theory and representation theory.
The representations of (affine) Weyl groups are well understood, while for general Coxeter groups they are still to be explored.

In the remarkable paper \cite{KL79}, Kazhdan and Lusztig introduced the notion of 2-sided cells for arbitrary Coxeter system $(W,S)$.
Each 2-sided cell gives a representation of $W$, called a cell representation.
When Lusztig's $\af$-function has a uniform bound on $W$, it is well known that any irreducible representation $V$ corresponds to a unique 2-sided cell $\mathcal{C}$ in the sense that $V$ is a quotient of the cell representation provided by $\mathcal{C}$.
We recollect these notions briefly in \S \ref{KL}.
Incidentally, we obtain a simple proof of the existence of the lowest 2-sided cell with the boundedness assumption on the function $\af$.

On the other hand, elements with $\af$-function value $1$ in an irreducible Coxeter system (i.e. a Coxeter system with connected Coxeter graph) form a 2-sided cell, denoted by $\mathcal{C}_1$.
One of the results in this article is a necessary and sufficient condition for a nontrivial irreducible representation $V$ to correspond to $\mathcal{C}_1$:
\begin{equation*}
 \begin{split}
  &\text{there is no $v \in V \setminus \{0\}$, such that $r \cdot v = t \cdot v = -v$} \\
  &\text{for some $r \neq t \in S$ with $m_{rt} < \infty$.}
 \end{split} \tag{A1}
\end{equation*}
See \ref{A1-condition} for details.
In particular, the geometric representation $\Vg$ satisfies \eqref{A1}.

There is an interesting class of representations ``between'' $\Vg$ and the ones satisfying \eqref{A1}, namely, the representations $V$ satisfying the following:
\begin{equation*}
 \begin{split}
   &\text{$V$ has a basis $\{\alpha_s \mid s \in S\}$, and for any $s \in S$ there is a subspace}  \\
   &\text{$H_s$ of codimension $1$ such that $s|_{H_s} = \operatorname{Id}_{H_s}$, and $s \cdot \alpha_s = -\alpha_s$.}
 \end{split} \tag{IR}
\end{equation*}
Let's refer to such representations as \emph{IR-representations}.
The main aim of \S \ref{secIR} is to present a classification of IR-representations assuming $m_{st} < \infty$ for any pair $s,t \in S$.
It turns out that IR-representations can be classified by some sets of natural numbers and one dimensional representations of the first homology groups of certain graphs.

In \S \ref{secR1}, we investigate the reducibility of IR-representations, which is encoded in a matrix defined by \eqref{matrixA}.
Such a criterion is similar to that of the geometric representation.
As an application, we determine all semisimple R-representations, which are defined to be a slightly larger class than IR-representations.
In \S \ref{sec-eg}, we will see that if $(W,S)$ is simply laced (i.e. $m_{st} = 2$ or $3$, $\forall s,t \in S$) and there is no more than one circuit in its Coxeter graph, all irreducible representations of $W$ corresponding to $\mathcal{C}_1$ are R-representations, thus we have found all irreducible quotients of the cell representation provided by $\mathcal{C}_1$ for such Coxeter groups.

Further we present some odds and ends in \S \ref{sec6}.
Proposition \ref{6.1} tells us when an IR-representation has a $W$-invariant bilinear (or sesquilinear) form, and \ref{6.4} discusses the dual of an IR-representation.

The story is similar when we are in the general setting.
In \S \ref{sec-gen}, we drop the assumption $m_{st} < \infty$ and present the classification of IR-representations.

In a sequel paper \cite{Hu22}, we use similar methods to construct some infinite dimensional irreducible representations of a large class of Coxeter groups, suggesting that the representation theory of general Coxeter groups is quite different from finite and affine counterparts.

\smallskip\noindent\textbf{Acknowledgement} The author would like to thank professor Nanhua Xi and Tao Gui for useful discussions.
The author is also grateful to professor Si'an Nie for pointing out an error of \ref{8.3} in a previous draft.

\section{Preliminaries on dihedral groups and graphs} \label{Dm-and-graph}
In the classification of IR-representations, we will investigate the way of ``gluing'' representations of dihedral groups, using certain graphs.
For convenience, let's recall in this section the representations of dihedral groups and integral homology of graphs.
The base field in this article is assumed to be $\mathbb{C}$ unless otherwise specified.
We use $e$ to denote the unity in a group.

\smallskip\noindent\textbf{Representations of finite dihedral groups} Let $D_m := \langle r,t \mid r^2 = t^2 = (rt)^m = e \rangle$ be a finite dihedral group, let $V_{r,t} := \mathbb{C}\beta_r \oplus \mathbb{C}\beta_t$.
For natural numbers $1 \leq k \leq m/2$, define the action $\rho_k : D_m \to \operatorname{GL}(V_{r,t})$ by
\begin{align*}
    r \cdot \beta_r & = - \beta_r, & r \cdot \beta_t & = \beta_t + 2\cos\frac{k\uppi}{m} \beta_r,\\
    t\cdot \beta_t & = - \beta_t, & t \cdot \beta_r &= \beta_r + 2\cos\frac{k\uppi}{m} \beta_t.
\end{align*}
Intuitively, $D_m$ acts on the plane via two reflections, the two reflection axis forming an angle of $k\uppi/m$, see Figure \ref{rho-k}.

\begin{figure}[ht]
    \centering
    \begin{tikzpicture}
      \draw[dashed] (-1.3,0)--(2,0);
      \draw[dashed] (240:1.3)--(60:2);
      \draw[->] (0,0)--(0,1.5);
      \draw[->] (0,0)--(330:1.5);
      \draw (0.3,0) arc (0:60:0.3);
      \draw[<->] ($(1.9,0) + (300:0.2)$) arc (300:420:0.2);
      \draw[<->] ($(60:1.9) + (0.2,0)$) arc (0:120:0.2);
      \node[left] (ar) at (0,1.5) {$\beta_r$};
      \node[below] (at) at (330:1.5) {$\beta_t$};
      \node[right] (r) at (2.1,0) {$r$};
      \node[right] (t) at (63:2.3) {$t$};
      \node[right] (angle) at (0.2,0.3) {$\frac{k\uppi}{m}$};
    \end{tikzpicture}
    \caption{$\rho_k : D_m \to \operatorname{GL}(V_{r,t})$}\label{rho-k}
\end{figure}

If $k < m/2$, $\rho_k$ is irreducible.
If $m$ is even and $k = m/2$, then $\rho_k = \varepsilon_r \oplus \varepsilon_t$ splits into a direct sum of two representations of dimension $1$, where
\begin{equation}\refstepcounter{Theorem} \label{eq-epsilon-rt}
    \varepsilon_r: r \mapsto -1, t \mapsto 1; \quad \quad
    \varepsilon_t: r \mapsto 1, t \mapsto -1.
\end{equation}
Denote by $\mathds{1}$ the trivial representation, by $\varepsilon$ the sign representation, i.e.
\begin{equation}\refstepcounter{Theorem}\label{eq-epsilon}
  \varepsilon: r \mapsto -1, t \mapsto -1.
\end{equation}

\begin{Lemma} \label{lem-Dm-rep}
All irreducible representations of $D_m$ are
\begin{gather*}
    \{\mathds{1}, \varepsilon, \rho_1, \dots, \rho_{\frac{m-1}{2}}\}, \quad \text{if $m$ is odd}; \\
    \{\mathds{1}, \varepsilon, \varepsilon_r, \varepsilon_t, \rho_1, \dots, \rho_{\frac{m}{2}-1}\}, \quad \text{if $m$ is even}.
\end{gather*}
\end{Lemma}

\begin{Remark} \label{rmk2.2}
The $+1$-eigenspaces of $r$ and $t$ in $\rho_k$ are both one dimensional.
However, any nonzero vector cannot be fixed by $r$ and $t$ simultaneously.
\end{Remark}

\begin{Remark}\label{rmk2.3}
By Schur's lemma, when $k < m/2$, an endomorphism of $\rho_k$ must be a scalar. This will be frequently used in the following sections.
\end{Remark}

\begin{Remark}\label{rmk2.4}
If $d>1$ is a common divisor of $k$ and $m$, then $\rho_k$ factors through $D_{\frac{m}{d}}$:
\begin{equation*}
    \xymatrix{D_m \ar[rr]^{\rho_k} \ar@{->>}[rd]& & \operatorname{GL}(V_{r,t}) \\
             & D_{\frac{m}{d}} \ar[ur]_{\rho_{\frac{k}{d}}} &}
\end{equation*}
\end{Remark}

\begin{Lemma} \label{lem-2.7}
When $1 \leq k < m/2$, there exists a bilinear form $\operatorname{B}$ on $V_{r,t}$, unique up to a $\mathbb{C}$-scalar, invariant under the action of $D_m$ via $\rho_k$:
\begin{equation}\refstepcounter{Theorem} \label{eq-Bform}
    \operatorname{B}(\beta_r,\beta_r) = \operatorname{B}(\beta_t, \beta_t) = 1, \quad \operatorname{B}(\beta_r, \beta_t) = \operatorname{B}(\beta_t, \beta_r) = -\cos \frac{k \uppi}{m}.
\end{equation}
In particular, $\operatorname{B}$ is symmetric.
\end{Lemma}

Let $\operatorname{H}: V_{r,t} \times V_{r,t} \to \mathbb{C}$ be the sesquilinear form on $V_{r,t}$ defined by \eqref{eq-Bform}, then it is invariant under the action $\rho_k$ of $D_m$. If $\operatorname{H}^\prime$ is another one invariant under $\rho_k$, then $\operatorname{H}^\prime$ is a $\mathbb{C}^\times$-scalar multiple of $\operatorname{H}$.
Further, $\operatorname{H}^\prime$ is Hermitian if and only if it is a $\mathbb{R}^\times$-scalar multiple of $\operatorname{H}$.

\smallskip\noindent\textbf{Representations of the Infinite dihedral group}
Let $D_\infty := \langle r,t \mid r^2 = t^2 = e \rangle$ be the infinite dihedral group.
There are four representations of dimension 1, namely, $\mathds{1}, \varepsilon, \varepsilon_r, \varepsilon_t$ (defined in the same way as in \eqref{eq-epsilon-rt}, \eqref{eq-epsilon}).

Let $\varrho_{x,y}: D_\infty \to \operatorname{GL}(V_{r,t})$ be the representation on $V_{r,t}$ defined by
\begin{align*}
    r \cdot \beta_r & = - \beta_r, & r \cdot \beta_t & = \beta_t + x \beta_r,\\
    t\cdot \beta_t & = - \beta_t, & t \cdot \beta_r &= \beta_r + y \beta_t.
\end{align*}
where $x, y \in \mathbb{C}$.
To distinguish from $\rho_k$'s above, we use the variation symbol $\varrho$.
\begin{Lemma}
  Suppose $x, y, x^\prime, y^\prime \in \mathbb{C}^\times$.
  \begin{enumerate}
    \item $\varrho_{x,y} \simeq \varrho_{x^\prime, y^\prime}$ if and only if $x y = x^\prime y^\prime$.
    \item $\varrho_{x,y}$ is irreducible if and only if $x y \neq 4$. 
  \end{enumerate}
\end{Lemma}

To simplify notation and eliminate ambiguity, for any $z \in \mathbb{C}$, choose and fix $u = u(z) \in \mathbb{C}$ such that $u^2 = z$, and if $z \in \mathbb{R}^+$ then $u$ is chosen to be positive. Write $\varrho_z := \varrho_{u,u}$.
When $z = 4 \cos^2 \frac{k \uppi}{m}$ for some $k,m \in \mathbb{N}, 1 \leq k < \frac{m}{2}$, $\varrho_z$ factors through $D_m$.
When $z = 4$, $\varrho_4 = \varrho_{2,2}$ is the geometric representation.
When $z = 0$, $\varrho_0 = \varrho_{0,0} \simeq \varepsilon_r \oplus \varepsilon_t$.

Besides, denote $\varrho_r^t := \varrho_{1,0}$ and $\varrho_t^r := \varrho_{0,1}$, then $\varrho_r^t, \varrho_t^r$ are both indecomposable.
$\varrho_r^t$ has a sub-representation $\varepsilon_r$ with quotient $\varepsilon_t$, while $\varrho_t^r$ has a sub-representation $\varepsilon_t$ with quotient $\varepsilon_r$.

There is an ``exotic'' representation (which we will not use in this article) of $D_\infty$ defined by $\varrho_{\varepsilon}^{\mathds{1}}: r \mapsto (\begin{smallmatrix}
  -1 & 0 \\
  0 & 1
  \end{smallmatrix}),
  t \mapsto (\begin{smallmatrix}
              -1 & 1 \\
              0 & 1
            \end{smallmatrix})$.
It has $\varepsilon$ as its sub-representation with quotient $\mathds{1}$.
\begin{Lemma}
  Above are all indecomposable representations of $D_\infty$ of dimension $1$ or $2$ (up to isomorphism), listed as
  \begin{equation*}
    \{\mathds{1}, \varepsilon, \varepsilon_r, \varepsilon_t, \varrho_r^t, \varrho_t^r, \varrho_{\varepsilon}^{\mathds{1}}\} \cup \{\varrho_z \mid z \in \mathbb{C}^\times\}.
  \end{equation*}
  All irreducible representations of $D_m$ are included, i.e.
  \begin{equation*}
    \{\mathds{1}, \varepsilon, \varepsilon_r, \varepsilon_t\} \cup \{\varrho_z \mid z \in \mathbb{C}^\times \setminus \{4\}\}.
  \end{equation*}
\end{Lemma}

\begin{Lemma}\label{lem2.11}
  Let $\varrho$ be an indecomposable representation of $D_\infty$ of dimension 2.
  \begin{enumerate}[resume]
    \item an endomorphism of $\varrho$ must be a scalar;
    \item \label{2.4.4} there is a $D_\infty$-invariant bilinear form $\operatorname{B}$ on $\varrho$, unique up to a scalar. The representing matrix of $\operatorname{B}$ (with respect to the basis $(\beta_r, \beta_t)$) is
        \begin{equation*}
          \begin{pmatrix}
            1 & -\frac{u}{2} \\
            -\frac{u}{2} & 1
          \end{pmatrix}, \text{if } \varrho = \varrho_z, z \neq 0;
          \begin{pmatrix}
            1 & 0 \\
            0 & 0
          \end{pmatrix}, \text{if } \varrho = \varrho_t^r;
          \begin{pmatrix}
            0 & 0 \\
            0 & 1
          \end{pmatrix}, \text{if } \varrho = \varrho_r^t \text{ or } \varrho_\varepsilon^\mathds{1}.
        \end{equation*}
  \end{enumerate}
\end{Lemma}

\smallskip\noindent\textbf{Graphs and their homology}
By definition, a \emph{(undirected) graph} $G = (S,E)$ consists of a set $S$ of vertices and a set $E$ of edges, elements in $E$ being of the form $\{s,t\}$ (unordered), where $s,t \in S$.

For our purpose, we only consider finite graphs without loops and multiple edges, i.e. there are no edges of the form $\{s,s\}$, each pair $\{s,t\}$ occurs at most once in $E$, and $S$ is a finite set.
Such a graph can be viewed as a finite simplicial complex in a natural way: vertices as $0$-simplices, edges as $1$-simplices.
A \emph{connected component} of the graph $G$ is defined to be that of this simplicial complex.

The chain complex (with coefficient $\mathbb{Z}$) of $G$ looks like:
\begin{equation*}
    0 \to \operatorname{C}_1(G) \xrightarrow{\partial} \operatorname{C}_0(G) \to 0.
\end{equation*}
The homology group $\operatorname{H}_1(G) = \ker \partial$ is a finitely generated free abelian group since $\operatorname{C}_1(G)$ is such a group.

Call $(s_1, s_2, \dots, s_n)$ a \emph{path} in $G$ if $\{s_i\}_i \subseteq S$ and $\{s_i, s_{i+1}\} \in E, \forall i$.
If $s_1 = s_n$, say such a path is a \emph{closed path}.
We don't distinguish the closed paths $(s_1, \dots, s_{n-1}, s_1)$ and $(s_2, \dots, s_{n-1}, s_1, s_2)$.
Further, if $s_1, \dots, s_{n-1}$ are distinct in this closed path, call it a \emph{circuit}.

For a connected graph $G$, a \emph{spanning tree} is defined to be a subgraph $T = (S,E_0)$ with the same vertex set $S$ and an edge set $E_0 \subseteq E$, such that $T$ is connected, and there are no circuits in $T$.
For general $G$, choosing a spanning tree for each component of $G$, their union is called a \emph{spanning forest} of $G$.
Spanning forests always exist, but are not unique in general.

Fix a spanning forest $G_0 = (S,E_0)$ of $G$.
For any edge $\mathfrak{e} \in E \setminus E_0$, fixing an orientation of $\mathfrak{e}$, there is a unique circuit in $(S, E_0 \cup \{\mathfrak{e}\})$ compatible with the orientation of $\mathfrak{e}$, denoted by $c_\mathfrak{e}$, representing an element in $\operatorname{H}_1(G)$.

\begin{Lemma}\label{lem2.9}
    $\operatorname{H}_1(G) = \bigoplus_{\mathfrak{e} \in E \setminus E_0} \mathbb{Z} c_\mathfrak{e}$.
\end{Lemma}

In a Coxeter graph $G$, the number $m_{st}$ is regarded as a label on the edge $\{s,t\}$, rather than a multiplicity.
Denote again by $\operatorname{H}_1(G)$ the first integral homology of $G$ forgetting labels on its edges.

\section{Backgrounds on Kazhdan-Lusztig cells} \label{KL}
This section explains partial motivations for considering IR-representations.
Fixing a Coxeter system $(W,S)$ of finite rank, let $\mathcal{H} := \mathcal{H}(W,S)$ be a free $\mathbb{C} [q^{\pm \frac{1}{2}}]$-module with basis $\{\nT_w\}_{w \in W}$.
Define a multiplication on $\mathcal{H}$ by
\begin{equation*}
  \nT_s \nT_w := \begin{cases}
                  \nT_{sw} , & \mbox{if } sw > w, \\
                  (\hq - \mhq) \nT_w + \nT_{sw} , & \mbox{if } sw < w,
                 \end{cases} \quad \forall s \in S, w \in W.
\end{equation*}
Then $\mathcal{H}$ forms an associative algebra with unity $\nT_e$, called the \emph{Hecke algebra} of $(W,S)$.
There is an involutive $\mathbb{C}$-algebra homomorphism on $\mathcal{H}$ defined by
\begin{equation*}
  \overline{\phantom{a}}: q^{\frac{1}{2}} \mapsto q^{-\frac{1}{2}}, \quad \nT_w \mapsto \nT_{w^{-1}}^{-1}.
\end{equation*}

\begin{Theorem}{\upshape (\cite[Theorem 1.1]{KL79})}
  There is a unique basis $\{C_w\}_{w \in W}$ of $\mathcal{H}$,
  \begin{equation*}
    C_w = \sum_{y \in W} (-1)^{\ell(w)+\ell(y)} q^{\frac{1}{2} (\ell(w) - \ell(y))} \overline{P_{y,w}} \nT_y,
  \end{equation*}
  such that $\forall y,w \in W$, $\overline{C_w} = C_w$, $P_{y,w} \in \mathbb{Z}[q]$, $P_{y,w} = 0$ unless $y \leq w$, $P_{w,w} = 1$, and if $y < w$, then $\deg_q P_{y,w} \leq \frac{1}{2} (\ell(w) - \ell(y) - 1)$.
\end{Theorem}

\begin{Examples} \label{eg3.4}
 \begin{enumerate}
   \item $C_e = \nT_e$; for $s \in S$, $C_s = \nT_s - \hq$.
   \item \label{3.4.3} If $r \neq t \in S$ with $m_{rt} = m < \infty$, write for simplicity $r_{k} := rtr \cdots$ (product of $k$ factors), $t_k = trt \cdots$ ($k$ factors), and $w_{rt} := r_{m} = t_{m}$. Then
       \begin{equation}\refstepcounter{Theorem} \label{Cwrt}
         C_{w_{rt}} = \nT_{w_{rt}} - \hq (\nT_{r_{m-1}} + \nT_{t_{m-1}}) + \dots + (-1)^i \hq[i] (\nT_{r_{m-i}} + \nT_{t_{m-i}}) + \dots + (-1)^m \hq[m].
       \end{equation}
   \item \label{3.4.1} $P_{y,w}(0) = 1$ for all $y \leq w$.
 \end{enumerate}
\end{Examples}

For any $x, y, w \in W$, let $h_{x, y, w} \in \AC$ such that $C_x C_y = \sum_{w \in W} h_{x,y,w} C_w$.
The following result is useful but highly nontrivial.

\begin{Lemma} \label{lemma-pos} {\upshape (\cite[Corollary 1.2]{EW14})}
  $\forall x,y,w \in W$, it holds $P_{y,w} \in \mathbb{N}[q]$. Write $h_{x,y,w} = \sum_{i \in \mathbb{Z}} c_i \hq[i]$, then $(-1)^i c_i \in \mathbb{N}$.
\end{Lemma}

For $x, y \in W$, say $y \lrleq x$ if $\exists H_1, H_2 \in \mathcal{H}$, such that $C_y$ has nonzero coefficient in the expression of $H_1 C_x H_2$ with respect to the basis $\{C_w\}_{w \in W}$.
Say $x \lrequiv y$ if $x \lrleq y$ and $y \lrleq x$.
It holds that $\lrleq$ is a pre-order on $W$, and $\lrequiv$ is an equivalence relation on $W$. The equivalence classes are called \emph{2-sided cells}.

Let $\af(w) := \min \{i \in \mathbb{N} \mid \hq[i] h_{x,y,w} \in \mathbb{Z}[\hq], \forall x, y \in W\}$, then $\af(w)$ is a well-defined natural number.
Say the function $\af$ is \emph{bounded} if there is $N \in \mathbb{N}$ such that $\af(w) \leq N$ for any $w \in W$.

\begin{Proposition} \label{prop3.8} {\upshape (\cite{Lusztig85-cell-i,Lusztig87-cell-ii, Lusztig14-hecke-unequal})}
  \begin{enumerate}
    \item \label{3.8.1} Let $s \in S$, $w \in W$, suppose $sw > w$, then $C_sC_w  = C_{sw} + \sum \mu_{y,w} C_y$ where the sum runs over all $y < w$ such that $\deg_q P_{y,w} = \frac{1}{2} (\ell(w) - \ell(y) - 1)$ and $s y < y$, and $\mu_{y,w}$ is the coefficient of the top-degree term in $P_{y,w}$.
        The formula is similar for the case $ws > w$.
        In particular, $sw \lrleq w$ if $sw > w$, $ws \lrleq w$ if $ws > w$.
    \item \label{3.8.3} $\af(w) = 0$ if and only if $w = e$. For any $s \in S$, $\af(s) = 1$.
    \item \label{3.8.5} If $m_{rt} < \infty$ for some $r \neq t \in S$, let $w = rtr\cdots$ ($m_{rt}$ factors), then $\af(w) = m_{rt}$.
    \item \label{3.8.4} If $y \lrleq x$, then $\af(y) \geq \af(x)$. Thus $y \lrequiv x$ yields $\af(y) = \af(x)$, and $\af(\mathcal{C})$ makes sense for any 2-sided cell $\mathcal{C}$.
    \item \label{3.8.6} Suppose $\af$ is bounded. If $\af(y) = \af(x)$ and $y \lrleq x$, then $y \lrequiv x$.
  \end{enumerate}
\end{Proposition}

Let $\mathcal{C}$ be a 2-sided cell, define $\mathcal{J}_\mathcal{C}$ to be a free $\AC$-module with a basis $\{J_w \mid w \in \mathcal{C}\}$.
The following formula defines an $\mathcal{H}$-module structure on $\mathcal{J}_\mathcal{C}$,
\begin{equation*}
  C_x \cdot J_y := \sum_{w \in \mathcal{C}} h_{x,y,w} J_w, \quad \forall x \in W, y \in \mathcal{C}.
\end{equation*}

By convention, a \emph{representation} of $\mathcal{H}$ indicates a $\mathbb{C}$ vector space with an $\mathcal{H}$-module structure, where $\hq$ acts by some $\mathbb{C}^\times$-scalar.
The following proposition is due to Lusztig.
We present a proof here briefly.

\begin{Proposition} \label{prop3.10}
  Suppose the function $\af$ is bounded, and $V \neq 0$ is an irreducible representation of $\mathcal{H}$, then there is a unique 2-sided cell $\mathcal{C}$, such that
  \begin{enumerate}
    \item for any $w \in W$, if $C_w \cdot V \neq 0$, then it holds $\mathcal{C} \lrleq w$;
    \item there is an element $w \in \mathcal{C}$ such that $C_w \cdot V \neq 0$.
  \end{enumerate}
  In other words, $V$ is a simple quotient of $\mathcal{J}_\mathcal{C}$.
\end{Proposition}

\begin{Proof}
  By \ref{prop3.8}\ref{3.8.6}, there exists a 2-sided cell $\mathcal{C}$ satisfying
  \begin{equation}\refstepcounter{Theorem} \label{P}
    \text{$\exists w \in \mathcal{C}$, such that $C_w \cdot V \neq 0$; $\forall y \lrleq \mathcal{C}$ and $y \notin \mathcal{C}$, $C_y \cdot V = 0$.}
  \end{equation}
  Take $v \in V$, such that $C_w \cdot v \neq 0$ for some $w \in \mathcal{C}$, then the map $\varphi: \mathcal{J}_\mathcal{C} \to V$, $J_y \mapsto C_y \cdot v$ is a well-defined surjective homomorphism of $\mathcal{H}$-modules.
  Thus $V \simeq \mathcal{J}_\mathcal{C} / \ker \varphi$.
  Suppose $\mathcal{C}^\prime$ is another 2-sided cell, and $\exists x \in \mathcal{C}^\prime$ such that $C_x \cdot V \neq 0$, then there exists $y \in \mathcal{C}$, such that $C_x \cdot J_y \neq 0$, i.e. there is $z \in \mathcal{C}$ such that $h_{x,y,z} \neq 0$, thus $z \lrleq x$, $\mathcal{C} \lrleq \mathcal{C}^\prime$.
  In particular, $\mathcal{C}^\prime$ doesn't satisfy \eqref{P}.
\end{Proof}

In the proof, the boundedness of $\af$ is only used to find $\mathcal{C}$ satisfying \eqref{P}.
Without this assumption, if such a 2-sided cell exists, then it is unique.

Immediately we have the following corollary on the lowest 2-sided cell, which was first studied for affine Weyl groups in \cite{Lusztig85-cell-i, Shi87}, and later for general Coxeter systems in \cite{Xi12} with a slightly stronger assumption.

\begin{Corollary}
  Suppose $N = \max \{\af(w) \mid w \in W\} <\infty$, then $\mathcal{C}_{\text{\upshape low}} := \{w \in W \mid \af(w) = N\}$ is a 2-sided cell, and for any $w \in W$ we have $\mathcal{C}_{\text{\upshape low}} \lrleq w$.
\end{Corollary}

\begin{Proof}
  Let $\varepsilon: s \mapsto -1, \forall s \in S$ be the sign representation, then for any $w \in W$, $\varepsilon(C_w) = (-1)^{\ell(w)} \sum_{y} P_{y,w}(1)$.
  By \ref{eg3.4}\ref{3.4.1} and \ref{lemma-pos}, this number is nonzero.
  Let $\mathcal{C}$ be the 2-sided cell attached to $\varepsilon$, then \ref{prop3.10} shows $\mathcal{C} \lrleq w$ for any $w$.
  By \ref{prop3.8}\ref{3.8.4}, we have $\mathcal{C} \subseteq \mathcal{C}_{\text{low}}$.
  The converse inclusion is deduced from \ref{prop3.8}\ref{3.8.6}.
\end{Proof}

We shall need the following lemma due to Matsumoto and Tits.

\begin{Lemma} \label{lemma3.12} {\upshape (\cite[Theorem 1.9]{Lusztig14-hecke-unequal})}
  Let $s_1 \cdots s_n$ and $s_1^\prime \cdots s_n^\prime$ be two reduced expressions of $w$, then we can obtain one from the other by finite steps of replacement of the form
    \begin{equation*}
          \underbrace{rtr\cdots} _{m_{rt} \text{ factors}} = \underbrace{trt\cdots} _{m_{rt} \text{ factors}} \quad (m_{rt} < \infty).
    \end{equation*}
\end{Lemma}

For irreducible $(W,S)$, there is a 2-sided cell with a simple description, formulated by Lusztig as following.

\begin{Lemma} \label{lem-C1}
  \begin{enumerate}
    \item For $w \in W \setminus \{e\}$, $\af(w) = 1$ if and only if $w$ has a unique reduced expression.
    \item Let $\mathcal{C}_1 := \{w \in W \mid \af(w) = 1\}$. If $(W,S)$ is irreducible, then $\mathcal{C}_1$ is a 2-sided cell. Moreover, for any $x \in W \setminus \{e\}$, we have $x \lrleq \mathcal{C}_1$.
  \end{enumerate}
\end{Lemma}

\begin{Proof}
  If $w$ has different reduced expressions, we deduce that $\af(w) \geq 2$ from \ref{lemma3.12} and \ref{prop3.8}\ref{3.8.1}\ref{3.8.5}\ref{3.8.4}.
  Conversely, if the reduced expression of $w$ is unique, suppose $sw < w$, \cite[Proposition 3.8]{Lusztig83-intrep} proved that $w \lrequiv s$, thus $\af(w) = 1$ by \ref{prop3.8}\ref{3.8.3}\ref{3.8.4}.

  The first part of (2) is proved in \cite[Proposition 3.8]{Lusztig83-intrep}. The second part is deduced from \ref{prop3.8}\ref{3.8.1}.
\end{Proof}

In general, decompose the Coxeter graph into connected components, say $S = \sqcup_i S_i$, then $\mathcal{C}_1 = \sqcup_i \mathcal{C}_{1,i}$ is a union of 2-sided cells, where $\mathcal{C}_{1,i} = \{w \in W \mid \af(w) = 1, \text{ and $w$ is a product of elements in $S_i$}\}$.

In the remainder of this section, we specialize $\hq$ to 1, then the Hecke algebra $\mathcal{H}$ becomes the group algebra $\mathbb{C}[W]$, and $\mathcal{J}_\mathcal{C}$ becomes a complex representation of $W$.
Denote the image of $C_w \in \mathcal{H}$ in $\mathbb{C}[W]$ by $C_w$ again.

\begin{Example}\label{Dm-cell}
  Use notations in \S \ref{Dm-and-graph}. $D_m$ has three 2-sided cells, $\{e\}$, $\{w_{rt}\}$ and $\mathcal{C}_1$.
  For $1 \leq k < m/2$, $\rho_k$ is an irreducible representation of dimension 2.
  On this plane, $D_m$ consists of $m$ rotations and $m$ reflections (regard $e$ as the rotation by zero angle).
  By \eqref{Cwrt}, $\rho_k(C_{w_{rt}})$ is the difference of the $m$ rotations and the $m$ reflections (up to a sign, depends on the parity of $m$).
  But the sum of the $m$ rotations is zero, so is the sum of the $m$ reflections, thus $\rho_k(C_{w_{rt}}) = 0$.
  One may further check that $C_{w_{rt}}$ acts by zero on other irreducible representations \emph{except} $\varepsilon$.
\end{Example}

The main result of this section is the following.

\begin{Theorem} \label{A1-condition}
Let $V$ be a representation of $W$, then $C_w$ acts by zero for all $w \in W$ with $\af(w) > 1$ if and only if $V$ satisfies the following condition,
\begin{equation*} \label{A1}
 \begin{split}
  &\text{there is no $v \in V \setminus \{0\}$, such that $r \cdot v = t \cdot v = -v$} \\
  &\text{for some $r \neq t \in S$ with $m_{rt} < \infty$.}
 \end{split} \tag{A1}
\end{equation*}
{\upshape Call such representations \emph{A1-representations}.}
\end{Theorem}

\begin{Proof}
  Suppose $r \neq t \in S$ with $m_{rt} = m < \infty$, write $w_{rt} = rtr \cdots$ ($m$ factors).
  Restrict $V$ to the subgroup $D_m := \langle r,t \rangle$, then $V$ (either finite or infinite dimensional) decomposes into irreducible representations of $D_m$.
  If $v \in V \setminus \{0\}$ spans a sign representation $\varepsilon$, then $C_{w_{rt}} \cdot v \neq 0$ as said in \ref{Dm-cell}.
  By \ref{prop3.8}\ref{3.8.5}, $\af(w_{rt}) = m > 1$.

  Conversely, suppose $V$ satisfies \eqref{A1}, then $C_{w_{rt}} \cdot V = 0$ for any pair $r,t$ as above.
  If $w \in W$, $\af(w) > 1$, by \ref{lemma3.12} and \ref{lem-C1}, $w = x w_{rt} y$ for some $r,t$ with $m_{rt} < \infty$, and $\ell(w) = \ell(x) + m_{rt} + \ell(y)$.
  Let's do induction on $\ell(w)$, the starting point is those elements of the form $w_{rt}$.
  If $x \neq e$, choose $s \in S$ such that $sx < x$, then $sw < w$ and $\af(sw) \ge 2$.
  We have by \ref{prop3.8}\ref{3.8.1} $C_w = C_s C_{sw} - \sum \mu_{z, sw} C_z$.
  In particular, those $z$'s occurred in the summation satisfy $z \lrleq sw$, thus $\af(z) \ge \af(sw) \geq 2$.
  Also, note that $\ell(z) < \ell(w)$, $\ell(sw) < \ell(w)$, by induction hypothesis, $C_{sw}$ and those $C_z$'s act on $V$ by zero, hence so is $C_w$.
  The case $y \neq e$ is similar.
\end{Proof}

In particular, if $(W,S)$ is irreducible, and $V$ is a nontrivial irreducible representation of $W$, then it satisfies \eqref{A1} if and only if $V$ is a quotient of $\mathcal{J}_{\mathcal{C}_1}$.

\section{Classification of IR-representations} \label{secIR}
Let $(W,S)$ be a Coxeter system of finite rank.
For simplicity, from here to \S \ref{sec6}, we assume:
\begin{equation}\refstepcounter{Theorem}\label{eq-4.1}
  m_{rt} < \infty, \quad \forall r,t \in S.
\end{equation}
A representation $V$ of $W$ is said to be an \emph{IR-representation} if:
\begin{equation*} \label{IR}
 \begin{split}
   &\text{$V$ has a basis $\{\alpha_s \mid s \in S\}$, and for any $s \in S$ there is a subspace}  \\
   &\text{$H_s$ of codimension $1$ such that $s|_{H_s} = \operatorname{Id}_{H_s}$, and $s \cdot \alpha_s = -\alpha_s$.}
 \end{split} \tag{IR}
\end{equation*}
The two letters represent ``independent reflections''.
Clearly the geometric representation satisfies \eqref{IR}, and IR-representations satisfy the condition \eqref{A1}.

\begin{Lemma}\label{lem4.2}
  Let $V$ be an IR-representation. For any pair $r \neq t \in S$,
  \begin{enumerate}
    \item \label{4.2-1} $\alpha_r, \alpha_t$ span a sub-representation of the dihedral subgroup $\langle r,t \rangle$;
    \item \label{4.2-3} As a representation of the dihedral group $\langle r,t \rangle$, the subspace $\langle \alpha_r, \alpha_t \rangle$ spanned by $\alpha_r, \alpha_t$ is isomorphic to $\rho_{k_{rt}}$ for some $1 \leq k_{rt} \leq m_{rt}/2$.
  \end{enumerate}
\end{Lemma}

\begin{Proof}
For each $s \in S$, $V = H_s \oplus \mathbb{C} \alpha_s$.
For any $v \in V$, write $v = v_s + a \alpha_s$, where $v_s \in H_s$, $a \in \mathbb{C}$, then $s \cdot v = v - 2a \alpha_s$.
\ref{4.2-1} is clear.
\ref{4.2-3} follows from \ref{lem-Dm-rep} and the assumption \eqref{IR}.
\end{Proof}

In view of \ref{lem4.2}\ref{4.2-3}, there exist $a_{r}^t, a_{t}^{r} \in \mathbb{C}^\times$ (but not unique), such that the isomorphism from $\rho_{k_{rt}}$ to $\langle \alpha_r, \alpha_t \rangle$ sends $\beta_r$ to $a_r^t \alpha_r$ and $\beta_t$ to $a_t^r \alpha_t$, i.e.
\begin{equation*}
    r \cdot \alpha_t = \alpha_t + 2 \frac{a_r^t}{a_t^r} \cos \frac{k_{rt}\uppi}{m_{rt}} \alpha_r, \quad \quad
    t\cdot \alpha_r = \alpha_r + 2 \frac{a_t^r}{a_r^t} \cos \frac{k_{rt}\uppi}{m_{rt}} \alpha_t.
\end{equation*}
Thus, the IR-representation is determined by the datum $(k_{rt}, a_r^t)_{r,t \in S, r \neq t}$.

Conversely, given the datum $(k_{rt}, a_r^t)_{r,t \in S, r \neq t}$ such that
\begin{equation*}
  k_{rt} \in \mathbb{N}, \quad 1 \leq k_{rt} = k_{tr} \leq \frac{m_{rt}}{2}, \quad a_r^t \in \mathbb{C}^\times, \quad \forall r,t \in S, r \ne t,
\end{equation*}
define for each $r \in S$ a linear transformation on the space $V := \bigoplus_{s \in S} \mathbb{C}\alpha_s$ by:
\begin{equation}\refstepcounter{Theorem} \label{eq-4.2}
    r \cdot \alpha_r = - \alpha_r, \quad r \cdot \alpha_t = \alpha_t + 2 \frac{a_r^t}{a_t^r} \cos \frac{k_{rt}\uppi}{m_{rt}} \alpha_r, \quad \forall t \neq r.
\end{equation}

\begin{Lemma} \label{lem4.4-defrep}
  Formula \eqref{eq-4.2} defines an IR-representation structure on $V$.
\end{Lemma}

\begin{Proof}
  First of all, $r^2$ acts by identity.
  Now consider the action of $(rt)^{m_{rt}}$.
  On $\langle \alpha_r, \alpha_t \rangle$,  $(rt)^{m_{rt}}$ acts by identity since the two dimensional subspace forms the representation $\rho_{k_{rt}}$.
  For $s \neq r,t$ (if exists), write $U := \langle \alpha_s, \alpha_r, \alpha_t \rangle$ be the $3$ dimensional subspace.
  By construction, there are vectors $v_1$ in $\langle \alpha_r, \alpha_s \rangle$ with nonzero $\alpha_s$-coefficient, and $v_2$ in $\langle \alpha_r, \alpha_t \rangle$ with nonzero $\alpha_t$-coefficient, both staying invariant under $r$.
  Thus $v_1$ and $v_2$ span a $2$ dimensional subspace of $U$.
  Similarly, there are $v_3$ and $v_4$ spanning a $2$ dimensional subspace of $U$ staying invariant under $t$.
  By dimension reason, there is a nonzero $v \in U$ fixed by $r$ and $t$ simultaneously.
  By \ref{rmk2.2}, $v \notin \mathbb{C} \alpha_r \oplus \mathbb{C} \alpha_t$.
  Thus, $\alpha_s = x_v v + x_r \alpha_r + x_t \alpha_t$ for some $x_v, x_r, x_t \in \mathbb{C}$.
  It follows $(rt)^{m_{rt}} \cdot \alpha_s = \alpha_s$.
  Now we see that $V$ is a representation of $W$ via the actions defined above.

  From the argument, we also see that for each $r \in S$ there are $|S|-1$ linearly independent vectors stay invariant under $r$, thus this representation satisfies \eqref{IR}.
\end{Proof}

We need to determine when two different data $(k_{rt}, a_r^t)_{r,t}$ and $(l_{rt}, b_r^t)_{r,t}$ define isomorphic representations.
Denote by $V_1 = \bigoplus_{s \in S} \mathbb{C} \alpha_s$ and $V_2 = \bigoplus_{s \in S} \mathbb{C} \alpha_s^\prime$ the two IR-representations defined by the two data respectively.
By \ref{lem4.2}, it obviously holds:

\begin{Lemma}
  $V_1 \simeq V_2$ only if $k_{rt} = l_{rt}$, $\forall r, t \in S$.
\end{Lemma}

Suppose $V_1, V_2$ are isomorphic via $\varphi: V_1 \xrightarrow{\sim} V_2$.
Define a new Coxeter system $(\widetilde{W}, S)$ with the same generator set $S$, and the order of $rt$ being $\widetilde{m}_{rt} := m_{rt} / d_{rt}$, where $d_{rt} := \gcd (m_{rt}, k_{rt})$.
Denote by $\widetilde{G} = (S, E)$ the Coxeter graph of $(\widetilde{W}, S)$.
We have a canonical homomorphism $W \twoheadrightarrow \widetilde{W}$.
By \ref{rmk2.4}, the representations $V_1$ and $V_2$ of $W$ factor through $\widetilde{W}$.
As an IR-representation of $\widetilde{W}$, $V_1$ is defined by the datum $(\widetilde{k}_{rt}, a_r^t)_{r,t}$, where $\widetilde{k}_{rt} := k_{rt} / d_{rt}$.
Similar for $V_2$.

The chain group $\operatorname{C}_1(\widetilde{G})$ is a free abelian group generated by oriented edges $(r,t)$.
The datum $(k_{rt}, a_r^t)_{r,t}$ gives a one dimensional representation of $\operatorname{C}_1(\widetilde{G})$ defined by
\begin{equation}\refstepcounter{Theorem}\label{eq-4.3}
    (r,t) \mapsto \frac{a_t^r}{a_r^t}.
\end{equation}
This restricts to the subgroup $\operatorname{H}_1 (\widetilde{G})$, denoted by $\chi_1$.
Similarly, $(r,t) \mapsto b_t^r/b_r^t$ defines another character of $\operatorname{H}_1 (\widetilde{G})$, denoted by $\chi_2$.

\begin{Lemma}\label{lem4.5}
  $\chi_1 =\chi_2$.
\end{Lemma}

\begin{Proof}
Recall that $\varphi: V_1 \to V_2$ is an isomorphism.
Let $(s_1, s_2, \dots, s_n, s_1)$ be a circuit in $\widetilde{G}$, write $\alpha_i := \alpha_{s_i}$ and $a_i^j := a_{s_i}^{s_j}$ for convenience, similar for $m_{ij}, \widetilde{m}_{ij}, k_{ij},$ $\alpha_i^\prime, b_i^j$, etc.
By assumption, $\varphi(\alpha_1) = x \alpha_1^\prime$ for some $x \in \mathbb{C}^\times$.
Since ${s_1, s_2}$ is an edge of $\widetilde{G}$, we have $k_{12} \neq m_{12}/2$, thus $\langle \alpha_1, \alpha_2 \rangle$ and $\langle \alpha_1^\prime,  \alpha_2^\prime \rangle$ form isomorphic irreducible representations $\rho_{k_{12}}$ of $\langle s_1, s_2 \rangle$.
By Schur's lemma, we have $\varphi(\alpha_2) = x (a_1^2 b_2^1 / a_2^1 b_1^2) \alpha_2^\prime$.
Do this recursively along the circuit, finally we get
\begin{equation*}
    \varphi(\alpha_1) = x \frac{a_1^2 a_2^3 \cdots a_{n-1}^n a_n^1 b_2^1 b_3^2 \cdots b_n^{n-1} b_1^n}{a_2^1 a_3^2 \cdots a_n^{n-1} a_1^n b_1^2 b_2^3 \cdots b_{n-1}^n b_n^1} \alpha_1^\prime.
\end{equation*}
Hence, compared with $\varphi(\alpha_1) = x \alpha_1^\prime$,
\begin{equation*}
    \frac{a_2^1 a_3^2 \cdots a_n^{n-1} a_1^n}{a_1^2 a_2^3 \cdots a_{n-1}^n a_n^1} = \frac{b_2^1 b_3^2 \cdots b_n^{n-1} b_1^n}{b_1^2 b_2^3 \cdots b_{n-1}^n b_n^1}.
\end{equation*}
The two sides are the images of the circuit under $\chi_1$ and $\chi_2$ respectively.
\end{Proof}

Fixing those numbers $(k_{rt})_{r,t}$, we have a well defined map $\Theta = \Theta((k_{rt})_{r,t})$ from the set of isomorphic classes
\begin{equation*}
  \{\text{IR-representations defined by } (k_{rt}, a_r^t)_{r,t} \mid a_r^t \in \mathbb{C}^\times, \forall r,t \in S\} / \simeq
\end{equation*}
to the set $\{$characters of $\operatorname{H}_1(\widetilde{G}) \}$ by \ref{lem4.5}.

\begin{Lemma}\label{lem4.7}
  $(k_{rt})_{r,t}$ fixed as above, $\Theta$ is surjective.
\end{Lemma}

\begin{Proof}
Fix a spanning forest $\widetilde{G}_0 = (S, E_0)$ of $\widetilde{G} = (S,E)$.
For any $\mathfrak{e} \in E \setminus E_0$, choose an orientation of $\mathfrak{e}$, let $c_\mathfrak{e}$ be defined as in \S \ref{Dm-and-graph}, then by \ref{lem2.9},
\begin{equation*}
    \operatorname{H}_1(\widetilde{G}) = \bigoplus_{\mathfrak{e} \in E \setminus E_0} \mathbb{Z} c_\mathfrak{e}.
\end{equation*}
Giving a character $\chi$ of $\operatorname{H}_1(\widetilde{G})$ is equivalent to giving numbers $\{ x_\mathfrak{e} \in \mathbb{C}^\times \mid \mathfrak{e} \in E \setminus E_0 \}$, and assign $\chi(c_\mathfrak{e}) = x_\mathfrak{e}$.

For any distinct $r,t \in S$, let
\begin{equation*}
    a_r^t = \begin{cases}
              x_{(t,r)}, & \mbox{if } \{t,r\} \in E \setminus E_0, \mbox{ and is oriented as } (t,r), \\
              1, & \mbox{otherwise.} \\

            \end{cases}
\end{equation*}
Then, under $\Theta$, the datum $(k_{rt}, a_r^t)_{r,t}$ gives the character $\chi$.
\end{Proof}

\begin{Lemma}\label{lem4.8}
  $(k_{rt})_{r,t}$ fixed as above, $\Theta$ is injective.
\end{Lemma}

\begin{Proof}
Suppose $\Theta (V_1) = \Theta (V_2)$, where $V_1 = \bigoplus_{s \in S} \mathbb{C} \alpha_s$, $V_2 = \bigoplus_{s \in S} \mathbb{C} \alpha_s^\prime$ are IR-representations defined by two data $(k_{rt}, a_r^t)_{r,t}$ and $(k_{rt}, b_r^t)_{r,t}$ respectively.
In arbitrary component of $\widetilde{G}$, choose a vertex $s_1$.
For any other vertex $s$ in the same component, choose a path $(s_1, s_2, \dots, s_n = s)$ connecting $s_1$ and $s$.
Use notations $\alpha_i, a_i^j$, etc as in the proof of \ref{lem4.5}.
Define a linear map $\varphi: V_1 \to V_2$ by
\begin{align*}
    \alpha_1 & \mapsto \alpha_1^\prime, \\
    \alpha_n & \mapsto \frac{a_1^2 a_2^3 \cdots a_{n-1}^n b_2^1 b_3^2 \cdots b_n^{n-1}}{a_2^1 a_3^2 \cdots a_n^{n-1} b_1^2 b_2^3 \cdots b_{n-1}^n} \alpha_n^\prime.
\end{align*}
We need to show $\varphi$ is well-defined, i.e. $\varphi(\alpha_n)$ is independent of the choice of the path (but $\varphi$ does depend on the choice of $s_1$ in each component).
For another path connecting $s_1$ and $s_n = s$, say $(s_1, s_p, s_{p-1}, \dots, s_{n+1}, s_n)$, $p \ge n$, the two paths form a closed path, i.e.  $(s_1, \dots s_n, s_{n+1}, \dots, s_p, s_1)$.
Since $\Theta(V_1) = \Theta(V_2)$, we have
\begin{equation*}
    \frac{a_2^1 \cdots a_n^{n-1} a_{n+1}^n \cdots a_p^{p-1} a_1^p}{a_1^2 \cdots a_{n-1}^n a_{n}^{n+1} \cdots a_{p-1}^p a_p^1} = \frac{b_2^1 \cdots b_n^{n-1} b_{n+1}^n \cdots b_p^{p-1} b_1^p}{b_1^2 \cdots b_{n-1}^n b_{n}^{n+1} \cdots b_{p-1}^p b_p^1}.
\end{equation*}
Rewrite the equation,
\begin{equation*}
    \frac{a_1^2 a_2^3 \cdots a_{n-1}^n b_2^1 b_3^2 \cdots b_n^{n-1}}{a_2^1 a_3^2 \cdots a_n^{n-1} b_1^2 b_2^3 \cdots b_{n-1}^n} = \frac{a_1^p a_p^{p-1} \cdots a_{n+1}^n b_p^1 b_{p-1}^p \cdots  b_n^{n+1}}{a_p^1 a_{p-1}^p \cdots a_n^{n+1} b_1^p b_p^{p-1} \cdots  b_{n+1}^n}.
\end{equation*}
This indicates that $\varphi(\alpha_n)$ is independent of the choice of the path.

Clearly, $\varphi$ is a linear isomorphism.
It remains to verify that $\varphi$ is a homomorphism of representations.
It suffices to check $\varphi(s \cdot \alpha_t) = s \cdot \varphi(\alpha_t), \forall s,t \in S$.
If $s = t$, this is obvious by definition.
If $s$ and $t$ are distinct and not adjacent in $\widetilde{G}$, i.e. $\widetilde{m}_{st} = 2$, then $s \cdot \alpha_t = \alpha_t$ and $s \cdot \alpha_t^\prime = \alpha_t^\prime$, the equation also holds.
Now assume $\widetilde{m}_{st} \geq 3$, then $s$ and $t$ are in the same connected component of $\widetilde{G}$.
Suppose $\varphi (\alpha_s) = x \alpha_s^\prime$, $x \in \mathbb{C}^\times$, then
\begin{equation*}
    \varphi (\alpha_t) = x \frac{a_s^t b_t^s}{a_t^s b_s^t} \alpha_t^\prime.
\end{equation*}
It is then straightforward to verify $\varphi(s \cdot \alpha_t) = s \cdot \varphi(\alpha_t)$ by definition.
\end{Proof}

Combine the discussions above, we have proved:

\begin{Theorem}\label{thm-IR}
  Isomorphic classes of IR-representations of $(W,S)$ one-to-one correspond to the set of data (where $\widetilde{G}$ is determined by $(k_{rt})_{r,t}$)
  \begin{equation*}
    \Bigl\{\bigl((k_{rt})_{r \neq t \in S}, \chi\bigr) \Bigm\vert  1 \leq k_{rt} = k_{tr} \leq \frac{m_{rt}}{2}, \forall r, t; \phantom{a} \chi \text{ is a character of } \operatorname{H}_1(\widetilde{G})\Bigr\}.
  \end{equation*}
\end{Theorem}

\begin{Remark}
 \begin{enumerate}
   \item Take all $k_{rt}$ to be $1$, and $\chi$ to be the trivial character, then the IR-representation obtained is just the geometric representation $\Vg$.
   \item This classification of IR-representations also works over the base field $\mathbb{R}$.
   \item Ignoring the labels on edges, the graph $\widetilde{G}$ is a subgraph of $G$ (the Coxeter graph of $(W,S)$) with the same vertex set. An edge $\{r,t\}$ in $G$ disappears in $\widetilde{G}$ if and only if $m_{rt}$ is even, and $k_{rt}$ is chosen to be $m_{rt}/2$.
       It may happen that $G$ is connected but $\widetilde{G}$ has several components.
 \end{enumerate}
\end{Remark}

From the proof of \ref{lem4.8}, we have the following corollary.

\begin{Corollary}\label{cor4.11}
  Let $V$ be an IR-representation, $\widetilde{G}$ be the corresponding graph.
  Let $g$ be the number of connected components of $\widetilde{G}$, then $\operatorname{End}_W (V) \simeq \mathbb{C}^{\oplus g}$.
  In particular, if $g = 1$, any endomorphism of $V$ is a scalar multiplication.
\end{Corollary}

\section{Reducibility and R-representations} \label{secR1}
Recall that if $(W,S)$ is reducible, the geometric representation $\Vg$ decomposes canonically into a direct sum of pieces corresponding to components of the Coxeter graph $G$.
When $(W,S)$ is irreducible, the reducibility and structure of $\Vg$ can be described by the following proposition.

\begin{Proposition} {\upshape (\cite[Ch.V \S 4 no.7]{Bourbaki-Lie456})}
 Assume $(W,S)$ is irreducible.
 \begin{enumerate}
   \item There is a unique up to scalar $W$-invariant bilinear form $\operatorname{B}$ on $\Vg$.
   \item $\Vg$ is irreducible if and only if the bilinear form $\operatorname{B}$ is non-degenerate.
   \item When $\Vg$ is reducible, it has a maximal sub-representation $V_0$ with trivial $W$-action, and the quotient $\Vg / V_0$ is irreducible.
 \end{enumerate}
\end{Proposition}

In this section we will see that IR-representations admit a similar criterion. As an application, a larger class of representations are classified.

\smallskip\noindent\textbf{Reducibility of IR-representations}
Let $V$ be an IR-representation of $(W,S)$ defined by datum $(k_{rt}, a_r^t)_{r,t}$, and $\widetilde{G}$ be the Coxeter graph attached to the numbers $(k_{rt})_{r,t}$, with corresponding Coxeter group $\widetilde{W}$, see \S \ref{secIR}.
If $S = \sqcup_i S_i$ is the decomposition of $\widetilde{G}$ into connected components, let $W_i$ and $\widetilde{W}_i$ denote the subgroup of $W$ and $\widetilde{W}$ generated by $S_i$ respectively.
The following are obvious.

\begin{Lemma}\label{6.2}
  \begin{enumerate}
    \item $\widetilde{W} = \prod_{i} \widetilde{W}_i$.
    \item $V = \bigoplus_i V_i$ where $V_i = \bigoplus_{s \in S_i} \mathbb{C} \alpha_s$ is a sub-representation of $W$.
    \item Each $V_i$ is an IR-representation of the Coxeter system $(W_i, S_i)$, defined by the datum $(k_{rt}, a_r^t)_{r,t \in S_i, r\neq t}$, and factors through $\widetilde{W}_i$.
    \item For $j \neq i$, $W_j$ acts trivially on $V_i$.
  \end{enumerate}
\end{Lemma}

Thus, we can focus our study on each direct summand.

\begin{Proposition}\label{prop5.3}
Suppose $\widetilde{G}$ is connected, and $U \subsetneq V$ is a sub-representation.
\begin{enumerate}
  \item If $U \ne 0$, the action of $W$ on $U$ is trivial.
  \item If $V$ is reducible, then $V$ has a maximal sub-representation $V_0$ with trivial $W$-action, and the quotient $V/V_0$ is irreducible.
      $V$ is indecomposable.
\end{enumerate}
\end{Proposition}

\begin{Proof}
Take $v \in U$, suppose $s \cdot v \neq v$ for some $s \in S$, then $v - s \cdot v \in \mathbb{C}^\times \alpha_s$ by the proof of \ref{lem4.2}, thus $\alpha_s \in U$.
If $r \in S$ is adjacent to $s$ in $\widetilde{G}$, then $\alpha_r \in U$ since $\langle \alpha_s, \alpha_r \rangle$ forms an irreducible representation of $\langle s, r \rangle$.
Inductively, we obtain $U = V$ which is a contradiction.
(2) follows from (1).
\end{Proof}

Label elements in $S$ to be $S = \{s_1, s_2, \dots, s_n\}$, and use notations $\alpha_i, k_{ij}, a_i^j,$ etc as in the proof of \ref{lem4.5}.
Write $v = \sum_{j} x_j \alpha_j$, $x_j \in \mathbb{C}$, then
\begin{align*}
    s_i \cdot v & = - x_i \alpha_i + \sum_{j \neq i} x_j \Bigl(\alpha_j + 2 \frac{a_i^j}{a_j^i} \cos \frac{k_{ij} \uppi}{m_{ij}} \alpha_i \Bigr) \\
     & = \Bigl(-x_i + \sum_{j \neq i} 2 x_j \frac{a_i^j}{a_j^i} \cos \frac{k_{ij} \uppi}{m_{ij}}\Bigr) \alpha_i + \sum_{j \neq i} x_j \alpha_j.
\end{align*}
If $v$ is fixed by all $s_i$, then
\begin{equation}\refstepcounter{Theorem} \label{eq6.1}
    x_i - \sum_{j \neq i}  x_j \frac{a_i^j}{a_j^i} \cos \frac{k_{ij} \uppi}{m_{ij}} = 0, \quad \forall i = 1,\dots,n.
\end{equation}
$V$ is reducible if and only if there exists nonzero vector $v \in V$ fixed by each $s_i$, if and only if the set of equations \eqref{eq6.1} has a nonzero solution for variables $\{x_i\}_i$, if and only if the following $n \times n$ matrix is singular:
\begin{equation}\refstepcounter{Theorem} \label{matrixA}
    A := \begin{pmatrix}
           1 & - \frac{a_1^2}{a_2^1} \cos \frac{k_{12} \uppi}{m_{12}} & \dots \\
           - \frac{a_2^1}{a_1^2} \cos \frac{k_{12} \uppi}{m_{12}} & 1 &  \\
           \vdots &  & \ddots
         \end{pmatrix}
\end{equation}
(All diagonal elements are $1$, and for any $i \neq j$ the element at $(i,j)$-position is $- \frac{a_i^j}{a_j^i} \cos \frac{k_{ij} \uppi}{m_{ij}}$. Note that $A$ is not symmetric in general.)
The co-rank of $A$ equals to the dimension of the solution space, as well as the dimension of the maximal sub-representation of $V$.
To conclude:

\begin{Proposition}\label{prop6.3}
  Suppose $\widetilde{G}$ is connected.
  \begin{enumerate}
    \item \label{6.3-1} $V$ is indecomposable.
        $V$ is irreducible if and only if the matrix $A$ defined in \eqref{matrixA} is invertible.
    \item \label{6.3-2} If $V$ is reducible, and $V_0$ is the maximal sub-representation, then the quotient $V/V_0$ is irreducible of dimension $\operatorname{rank}(A)$.
  \end{enumerate}
  In general, $V$ is semisimple if and only if $A$ is invertible.
\end{Proposition}

\begin{Remark}
Fixing $(k_{rt})_{r,t}$ such that $\widetilde{G}$ is connected and adapt the choice in \ref{lem4.7}, the corresponding IR-representations are parameterized by a torus $(\mathbb{C}^\times)^N$, where $N$ is the rank of $\operatorname{H}_1 (\widetilde{G})$.
Those reducible representations form a subvariety of $(\mathbb{C}^\times)^N$.
\end{Remark}

\noindent\textbf{R-representations}
In the condition \eqref{IR}, we require that those ``reflection vectors'' $\{\alpha_s\}_s$ are linearly independent.
Now let's remove this restriction, but we still want to stay in the \eqref{A1} assumption, so we define \emph{R-representations} to be  representations $V$ of $W$ satisfying the following condition:
\begin{equation*}\label{R}
  \begin{split}
      & \text{$V$ is spanned by nonzero vectors $\{\alpha_s \mid s \in S\}$, and $\forall s \in S$, $\exists$ a subspace} \\
       & \text{$H_s$ of codimension $1$ such that $s|_{H_s} = \operatorname{Id}_{H_s}$, and $s \cdot \alpha_s = -\alpha_s$.} \\
       & \text{Moreover, $\forall r \neq t \in S$ with $m_{rt} < \infty$, $\alpha_r$ and $\alpha_t$ are linearly independent.}
  \end{split}\tag{R}
\end{equation*}
Note that at present we require $m_{rt} < \infty$ for any $r,t$ (see the beginning of \S \ref{secIR}), so each pair $\alpha_r, \alpha_t$ are not proportional.

By definition, R-representations form a larger class than IR-representations.
But in fact they don't produce too much more.

\begin{Theorem}\label{thm-R1}
Let $V$ be an R-representation.
\begin{enumerate}
    \item \label{5.6.1} There is a unique (up to isomorphism) IR-representation $V^\prime$ such that $V$ is a quotient of $V^\prime$, say, $\pi : V^\prime \twoheadrightarrow V$, and $W$ acts on $\ker \pi$ trivially.
    \item Isomorphism classes of semisimple R-representations one-to-one correspond to isomorphism classes of IR-representations, thus to the set of data $\{((k_{rt})_{r,t}, \chi)\}$ in \ref{thm-IR}.
        In particular, simple R-representations correspond to those data such that $\widetilde{G}$ is connected.
  \end{enumerate}
\end{Theorem}

\begin{Proof}
For such $V$, $V = \sum_{s \in S} \mathbb{C} \alpha_s$, and $\forall s \in S$, $V = H_s \oplus \mathbb{C} \alpha_s$.
By assumption \eqref{eq-4.1}, for any $r \neq t \in S$, $\langle \alpha_r, \alpha_t \rangle$ is two dimensional, and the same argument in \ref{lem4.2} still works.
Thus, we can extract the datum $(k_{rt}, a_r^t)_{r,t}$ from $V$ as in \S \ref{secIR}.
Let $V^\prime = \bigoplus_{s \in S} \mathbb{C} \alpha_s^\prime$ be the IR-representation defined by this datum, and $\pi: V^\prime \to V$ be the surjective linear map defined by $\alpha_s^\prime \mapsto \alpha_s$.
One may easily verify that $\pi$ is a homomorphism of representations, thus $V$ is a quotient of $V^\prime$.

Suppose $V^{\prime \prime} = \bigoplus_{s \in S} \mathbb{C} \alpha_s^{\prime \prime}$ is another IR-representation defined by $(l_{rt}, b_r^t)_{r,t}$, and $\pi^{\prime \prime}: V^{\prime \prime} \twoheadrightarrow V$ is a surjective homomorphism.
After rescaling, we may assume $\pi^{\prime \prime}(\alpha_s^{\prime \prime}) = \alpha_s$, $\forall s$.
Restrict to the dihedral subgroup $\langle r, t \rangle$, it must hold $k_{rt} = l_{rt}$, since $\langle \alpha_r^{\prime \prime}, \alpha_t^{\prime \prime} \rangle$ and $\langle \alpha_r, \alpha_t \rangle$ should be isomorphic as representations of $\langle r, t \rangle$.
Next, consider $\pi^{\prime \prime} (r \cdot \alpha_t^{\prime \prime})=r \cdot \alpha_t$, if $k_{rt} \neq m_{rt}/2$, then we obtain $a_r^t / a_t^r = b_r^t / b_t^r$ by \eqref{eq-4.2}.
Thus the two data define the same character of $\operatorname{H}_1 (\widetilde{G})$, and $V^\prime \simeq V^{\prime \prime}$ by \ref{thm-IR}, where $\widetilde{G}$ is the Coxeter graph attached to $(k_{rt})_{r,t}$. This proves the uniqueness part in (1).

Let $S = \sqcup_i S_i$ be the decomposition of $\widetilde{G}$ into components, then $V^\prime = \bigoplus_i V_i^\prime$ where $V_i^\prime = \bigoplus_{s \in S_i} \mathbb{C} \alpha_s^\prime$ is a sub-representation in $V^\prime$, as in \ref{6.2}.
Write $V_{i,0}^\prime$ to be the maximal sub-representation of $V_i^\prime$ (may be zero).
Note that $\ker \pi$ is a sub-representation of $V^\prime$ but $\ker \pi$ can not contain any $V_i^\prime$.
Suppose $0 \neq v \in \ker \pi$, write $v = \oplus_i v_i$ where each $v_i \in V_i^\prime$.
Take $s \in S_j$, for any $i \neq j$, we have $s \cdot v_i = v_i$.
If $s \cdot v_j \neq v_j$, the argument in \ref{prop5.3} shows $V_j^\prime \subseteq \ker \pi$, which is impossible.
Thus for each $i$, $v_i \in V_{i,0}^\prime$, $W$ acts on $\ker \pi$ trivially.
(1) is proved.

We have shown that $\ker \pi \subseteq V_0^\prime := \bigoplus_i V_{i,0}^\prime$.
On one hand, by \ref{prop6.3}, $V^\prime / V_0^\prime$ is a semisimple R-representation for any IR-representation $V^\prime$.
On the other hand, if $\ker \pi \subsetneq V_0^\prime$, suppose $V_{j,0}^\prime \nsubseteq \ker \pi$, then $\pi(V_j^\prime)$ is an indecomposable but reducible sub-representation of $V$.
Thus $V$ is not semisimple.
(2) follows.
\end{Proof}

\begin{Remark}\label{rmk5.7}
  \begin{enumerate}
    \item Let $v = \oplus_i v_i \in \ker \pi$ as in the proof, it may happen that $v_i \notin \ker \pi$.
        In fact, take any subspace $K$ of $\bigoplus_i V_{i,0}^\prime$, the quotient $V^\prime / K$ is an R-representation of $W$.
    \item \label{5.7.2} Let $\pi_1 : V^\prime \twoheadrightarrow V_1$, $\pi_2 : V^\prime \twoheadrightarrow V_2$, where $V^\prime$ is an IR-representation and $V_1, V_2$ are R-representations.
        Suppose $\widetilde{G}$ is connected, the same argument of \ref{cor4.11} shows $\operatorname{End}_W(V_1) = \mathbb{C}$ and $\operatorname{Hom}_W(V^\prime, V_1) \simeq \mathbb{C}$.
        It yields that $V_1 \simeq V_2$ if and only if $\ker \pi_1 = \ker \pi_2$.
        Thus for such $V^\prime$, the isomorphism classes of quotients of $V^\prime$ is parameterized by subspaces of the maximal sub-representation of $V^\prime$ (a union of Grassmannians).
        See \ref{7.2} for an example for the case $\widetilde{G}$ is disconnected.
  \end{enumerate}
\end{Remark}

\section{Additional properties} \label{sec6}
\textbf{$W$-invariant bilinear form}
First we determine when an IR-representation admits a $W$-invariant bilinear form.
By \ref{6.2}, we may focus on a direct summand of $V$.

\begin{Proposition}\label{6.1}
  Let $V$ be an IR-representation defined by datum $(k_{rt}, a_r^t)_{r,t}$, and $\chi$ be the corresponding character of $\operatorname{H}_1 (\widetilde{G})$. Assume $\widetilde{G}$ is connected.
  \begin{enumerate}
    \item $V$ has a nonzero $W$-invariant bilinear form if and only if $\chi (\operatorname{H}_1 (\widetilde{G})) \subseteq \{\pm 1\}$.
    \item $V$ has a nonzero $W$-invariant sesquilinear form if and only if $\chi (\operatorname{H}_1 (\widetilde{G})) \subseteq S^1$.
  \end{enumerate}
\end{Proposition}

\begin{Proof}
Suppose there is a nonzero $W$-invariant bilinear form $\operatorname{B}$ on $V$.
For any path $(s_1, s_2, \dots, s_n)$ in $\widetilde{G}$, suppose $\operatorname{B} (\alpha_1, \alpha_1) = x$ (here $\alpha_1 = \alpha_{s_1}$, similarly for other notations), then by \ref{lem-2.7} and the construction of $V$, it holds $\operatorname{B} (\alpha_2, \alpha_2) = x (a_1^2/a_2^1)^2$.
Do this recursively along this circuit, finally we get
\begin{equation}\refstepcounter{Theorem}\label{eq-6.1}
    \operatorname{B} (\alpha_n, \alpha_n) = x \left( \frac{a_1^2 a_2^3 \cdots a_{n-1}^n}{a_2^1 a_3^2 \cdots a_n^{n-1}} \right)^2.
\end{equation}
If $x = 0$, then $\operatorname{B}(\alpha_s, \alpha_s) = 0$ for any $s \in S$ since $\widetilde{G}$ is connected.
As a result, if $s,t \in S$ such that $\widetilde{m}_{st} \ge 3$, then $\operatorname{B}(\alpha_s, \alpha_t) = \operatorname{B}(\alpha_t, \alpha_s) = 0$ by \ref{lem-2.7}.
If $\widetilde{m}_{st} =2$, then $\operatorname{B}(\alpha_s, \alpha_t) = \operatorname{B}(s \cdot \alpha_s, s\cdot \alpha_t) = - \operatorname{B}(\alpha_s, \alpha_t)$ yields $\operatorname{B}(\alpha_s, \alpha_t) = 0$.
Therefore $\operatorname{B} = 0$ which is a contradiction. Thus $x \ne 0$.
Suppose now $s_1 = s_n$, then $\operatorname{B}(\alpha_n, \alpha_n) = x$, \eqref{eq-6.1} implies the image of the closed path under $\chi$ lies in $\{\pm 1\}$.

Conversely, suppose $\chi (\operatorname{H}_1 (\widetilde{G})) \subseteq \{\pm 1\}$.
Choose and fix a vertex $s_1$ in $\widetilde{G}$, define $\operatorname{B} (\alpha_1, \alpha_1) := 1$.
For another vertex $s$, choose a path $(s_1, s_2, \dots, s_n = s)$, in view of \ref{lem-2.7}, define
\begin{equation*}
    \operatorname{B} (\alpha_s, \alpha_s) := \left( \frac{a_1^2 a_2^3 \cdots a_{n-1}^n}{a_2^1 a_3^2 \cdots a_n^{n-1}} \right)^2.
\end{equation*}
For $t \neq s$, define
\begin{gather*}
    \operatorname{B} (\alpha_t, \alpha_s) = \operatorname{B} (\alpha_s, \alpha_t) := - \frac{a_s^t}{a_t^s} \cos \frac{k_{st} \uppi}{m_{st}} \operatorname{B} (\alpha_s, \alpha_s) \\
    \bigl(\text{also equals to } - \frac{a_t^s}{a_s^t} \cos \frac{k_{st} \uppi}{m_{st}} \operatorname{B} (\alpha_t, \alpha_t)\bigr).
\end{gather*}
$\operatorname{B}$ is well-defined, i.e. $\operatorname{B} (\alpha_s, \alpha_s)$ is independent of the choice of the path.
This is because, choosing another path $(s_1, s_p, s_{p-1}, \dots, s_{n+1}, s_n = s)$, $p \ge n$, the two paths form a circuit $(s_1, \dots, s_n, s_{n+1}, \dots, s_p, s_1)$, and our assumption says
\begin{equation*}
    \pm 1 = \frac{a_1^2 \cdots a_{n-1}^n a_n^{n+1} \cdots a_{p-1}^p a_p^1}{a_2^1 \cdots a_n^{n-1} a_{n+1}^n \cdots a_p^{p-1} a_1^p}, \text{ thus } \left( \frac{a_1^2 a_2^3 \cdots a_{n-1}^n}{a_2^1 a_3^2 \cdots a_n^{n-1}} \right)^2 = \left( \frac{a_1^p a_p^{p-1} \cdots a_{n+1}^n}{a_p^1 a_{p-1}^p \cdots a_n^{n+1}} \right)^2.
\end{equation*}
Now we need to check $\operatorname{B}$ is $W$-invariant. It suffices to check $\operatorname{B} (s \cdot \alpha_r, s \cdot \alpha_t) = \operatorname{B} (\alpha_r, \alpha_t)$ for $s, r, t \in S$.
If $s = r = t$, this is obvious.
If $s = r \neq t$, or $s = t \neq r$, or $s \neq r = t$, everything happens in a dihedral world, the equation holds by \ref{lem-2.7}.
If $s, r, t$ are distinct,
\begin{align*}
    \operatorname{B} (s \cdot \alpha_r, s \cdot \alpha_t) & = \operatorname{B} \Bigl(\alpha_r + 2 \frac{a_s^r}{a_r^s} \cos \frac{k_{sr} \uppi}{m_{sr}} \alpha_s, \alpha_t + 2 \frac{a_s^t}{a_t^s} \cos \frac{k_{st} \uppi}{m_{st}} \alpha_s\Bigr) \\
    & = \operatorname{B} (\alpha_r, \alpha_t) + 2 \frac{a_s^t}{a_t^s} \cos \frac{k_{st} \uppi}{m_{st}} \operatorname{B} (\alpha_r, \alpha_s) \\
    & \mathrel{\phantom{=}} + 2 \frac{a_s^r}{a_r^s} \cos \frac{k_{sr} \uppi}{m_{sr}} \operatorname{B} (\alpha_s, \alpha_t) + 4 \frac{a_s^t a_s^r}{a_t^s a_r^s} \cos \frac{k_{st} \uppi}{m_{st}} \cos \frac{k_{sr} \uppi}{m_{sr}} \operatorname{B} (\alpha_s, \alpha_s)\\
    & = \operatorname{B} (\alpha_r, \alpha_t).
\end{align*}

The proof of (2) is similar to the argument above.
\end{Proof}

\begin{Remark}
  \begin{enumerate}
    \item From the proof we see that the $W$-invariant bilinear form $\operatorname{B}$ (if exists) is symmetric and is unique up to a $\mathbb{C}^\times$-scalar.
    \item If we don't assume $\widetilde{G}$ is connected, write $V = \bigoplus_i V_i$ as in \ref{6.2}, then \ref{6.1} is applied on each $V_i$.
    \item By \ref{6.1}, there are only finitely many IR-representations of $W$ admitting an invariant bilinear form. The geometric representation is one of them.
  \end{enumerate}
\end{Remark}

Let $V$ be an IR-representation with connected $\widetilde{G}$.
Suppose $V$ is not irreducible, let $V_0$ be the maximal sub-representation.
Assume further $V$ admits a nonzero $W$-invariant bilinear form $\operatorname{B}$.
Let $v = \sum_{s} x_s \alpha_s \in V_0$, $x_s \in \mathbb{C}$, then
\begin{equation}\refstepcounter{Theorem}\label{eq-6.4}
    \operatorname{B} (v, \alpha_r) = x_r \operatorname{B} (\alpha_r, \alpha_r) - \sum_{s \neq r} x_s \frac{a_r^s}{a_s^r} \cos \frac{k_{sr} \uppi}{m_{sr}} \operatorname{B} (\alpha_r, \alpha_r) = 0, \quad \forall r \in S.
\end{equation}
The second equality is due to  \eqref{eq6.1}.
Let $V_1$ be an R-representation which is a quotient of $V$, $\pi : V \twoheadrightarrow V_1$ be the projection, then $\ker \pi \subseteq V_0$.
Formula \eqref{eq-6.4} tells us the bilinear form $\operatorname{B}$ descents to $V_1$.

Conversely, let $V_1$ be an R-representation with connected $\widetilde{G}$.
Suppose $V_1$ has a nonzero $W$-invariant bilinear form $\operatorname{B}$, by the same argument of \ref{6.1}, we shall see that $\chi (\operatorname{H}_1 (\widetilde{G})) \subseteq \{\pm 1\}$, and such bilinear form is unique up to a scalar.
Also, $\operatorname{B}$ can be lifted to $V$, the corresponding IR-representation.

\smallskip\noindent\textbf{Dual of an IR-representation}
Let $V = \bigoplus_{s \in S} \mathbb{C} \alpha_s$ be an IR-representation defined by datum $(k_{rt}, a_r^t)_{r,t}$, and $\{\alpha_s^*\}_{s \in S}$ be the dual basis of $\{\alpha_s\}_s$ in $V^* := \operatorname{Hom}_\mathbb{C} (V, \mathbb{C})$.
$V^*$ forms a representation of $W$ via $(s \cdot \alpha_t^*) (-) = \alpha_t^*( s \cdot -)$, $\forall s\in S$.
A direct computation shows
\begin{align}
    s \cdot \alpha_s^* & = - \alpha_s^* + \sum_{t \neq s} 2 \frac{a_s^t}{a_t^s} \cos \frac{k_{st} \uppi}{m_{st}} \alpha_t^*, \refstepcounter{Theorem}\label{eq6.5} \\
    s \cdot \alpha_t^* & = \alpha_t^*, \quad \forall t \neq s. \refstepcounter{Theorem}\label{eq6.6}
\end{align}
From this we can see that the unique (up to a scalar) $-1$-eigenvector of $s$ is
\begin{equation}\refstepcounter{Theorem}\label{eq-6.5}
    \gamma_s := \alpha_s^* - \sum_{t \neq s} \frac{a_s^t}{a_t^s} \cos \frac{k_{st} \uppi}{m_{st}} \alpha_t^*.
\end{equation}

\begin{Proposition} \label{6.4}
  \begin{enumerate}
    \item $\{\gamma_s\}_{s \in S}$ is linearly independent if and only if the matrix $A$ defined in \eqref{matrixA} is invertible, if and only if $V$ is semisimple.
    \item When $V$ is semisimple, $V^* = \bigoplus_{s} \mathbb{C} \gamma_s$ is an IR-representation of $W$.
        Suppose $((k_{rt})_{r,t},\chi)$ is the datum corresponding to  $V$, then $V^*$ corresponds to the datum $((k_{rt})_{r,t},\chi^*)$, here $\chi^*$ is the dual character, $\chi^* (-) = \chi (-)^{-1}$.
    \item In general case, $\{\gamma_s\}_{s \in S}$ spans a sub-representation $V_1 := \sum_{s \in S} \mathbb{C} \gamma_s$ of $V^*$.
        Moreover, $V_1$ is an R-representation isomorphic to the semisimple quotient of the IR-representation corresponding to the datum $((k_{rt})_{r,t}, \chi^*)$.
        The action of $W$ on the quotient $V^* / V_1$ is trivial.
  \end{enumerate}
\end{Proposition}

\begin{Proof}
By \eqref{eq-6.5}, the transition matrix from $\{\alpha_s^*\}_s$ to $\{\gamma_s\}_s$ is $A^\mathrm{T}$, the transpose of $A$, thus $\{\gamma_s\}_s$ is linearly independent if and only if $A$ is invertible.
The second part of (1) is due to \ref{prop6.3}.

When $\{\gamma_s\}_s$ is linearly independent, clearly $V^* = \bigoplus_{s} \mathbb{C} \gamma_s$, and $V^*$ is an IR-representation.
A direct computation shows
\begin{equation}\refstepcounter{Theorem}\label{eq-6.7}
    s \cdot \gamma_t = \gamma_t + 2 \frac{a_t^s}{a_s^t} \cos \frac{k_{st} \uppi}{m_{st}} \gamma_s,
\end{equation}
thus the datum $(k_{rt}, b_r^t)_{r,t}$ defines $V^*$, where $b_r^t = a_t^r$.
Therefore $V^*$ corresponds to the character $\chi^*$ of $\operatorname{H}_1(\widetilde{G})$.
(2) is proved.

In general, let $V^\prime$ be the IR-representation defined by $(k_{rt}, b_r^t)_{r,t}$, then $A^\mathrm{T}$ is the matrix of $V^\prime$ defined as in \eqref{matrixA}.
Notice that $V_1$ is an R-representation, and $\dim V_1$ equals to the rank of $A^\mathrm{T}$, thus $V_1$ is the semisimple quotient of $V^\prime$.
By \eqref{eq6.5} and \eqref{eq-6.5}, we have $s \cdot \alpha_s^* = \alpha_s^* -2 \gamma_s$.
Combined with \eqref{eq6.6}, we see that $W$ acts on $V^*/V_1$ trivially.
(3) is proved.
\end{Proof}

\section{General situation: allowing \texorpdfstring {$m_{st} = \infty$} {infinity}} \label{sec-gen}
From now on, we consider arbitrary Coxeter system $(W,S)$ of finite rank with Coxeter graph $G$, dropping the assumption that $m_{rt} < \infty$, $\forall r,t \in S$.
For any pair $r \neq t \in S$, define a parameter set
\begin{equation*}
    P_{rt} := P_{m_{rt}} = \begin{cases}
            \{\rho_k \mid k \in \mathbb{N}, 1 \leq k \leq \frac{m_{rt}}{2}\}, & \mbox{if } m_{rt} < \infty; \\
            \{\varrho_r^t, \varrho_t^r\} \cup \{\varrho_z \mid z \in \mathbb{C}\}, & \mbox{if } m_{rt} = \infty.
          \end{cases}
\end{equation*}
(The notations $\varrho$'s are defined in \S \ref{Dm-and-graph}.)
Note that $P_\infty$ consists of all IR-representations of $D_\infty$.
Suppose we are given $\delta_{rt} \in P_{rt}$ for any pair $r \neq t \in S$, such that $\delta_{rt} = \delta_{tr}$ representing the same representation of $\langle r,t \rangle$, define a Coxeter system $(\widetilde{W},S)$ with Coxeter graph $\widetilde{G}$ as follows.
\begin{enumerate}
  \item If $m_{rt} < \infty$ and $\delta_{rt} = \rho_{k}$, let $\widetilde{m}_{rt} = m_{rt} / d$ where $d = \operatorname{g.c.d.} (m_{rt}, k)$;
  \item If $m_{rt} = \infty$, $\delta_{rt} = \varrho_z$ and $z = 4 \cos^2 \frac{k \uppi}{m}$ for some co-prime $k,m \in \mathbb{N}$, $1 \leq k < m/2$, let $\widetilde{m}_{rt} = m$;
  \item If $m_{rt} = \infty$, $\delta_{rt} = \varrho_0$, let $\widetilde{m}_{rt} = 2$;
  \item Otherwise, let $\widetilde{m}_{rt} = \infty$.
\end{enumerate}

\begin{Theorem}
  Isomorphism classes of IR-representations of $(W,S)$ one-to-one correspond to the set of data
  \begin{equation*}
    \Bigl\{\bigl((\delta_{rt})_{r \neq t \in S}, \chi\bigr) \Bigm\vert \delta_{rt} = \delta_{tr} \in P_{rt}, \forall r, t; \phantom{a} \chi \text{ is a character of } \operatorname{H}_1(\widetilde{G})\Bigr\}.
  \end{equation*}
\end{Theorem}

\begin{Proof}
Let $V = \bigoplus_s \mathbb{C} \alpha_s$ be an IR-representation of $W$, \ref{lem4.2}\ref{4.2-1} still holds.
For any pair $r \neq t \in S$, there is an element $\delta_{rt} \in P_{rt}$, such that $\langle \alpha_r, \alpha_t \rangle$ forms a representation of $\langle r,t \rangle$ isomorphic to $\delta_{rt}$.
There exist $a_r^t, a_t^r \in \mathbb{C}^\times$ such that the isomorphism sends $a_r^t \alpha_r$ to $\beta_r$, $a_t^r \alpha_t$ to $\beta_t$.
The structure of $V$ is determined by the datum $(\delta_{rt}, a_r^t)_{r,t}$, since the action $r \cdot \alpha_t$ has been clear, for any $r,t$.

Conversely, given datum $(\delta_{rt}, a_r^t)_{r,t}$ in which $\delta_{rt} = \delta_{tr} \in P_{rt}$, $a_r^t \in \mathbb{C}^\times$, identify $\langle \alpha_r, \alpha_t \rangle$ and $\delta_{rt}$ via $a_r^t \alpha_r \mapsto \beta_r$, $a_t^r \alpha_t \mapsto \beta_t$, then $\langle r,t \rangle$ acts on $\langle \alpha_r, \alpha_t \rangle$.
The same argument in \ref{lem4.4-defrep} shows that this defines a $W$-action on $V := \bigoplus_s \mathbb{C} \alpha_s$.

For a datum $(\delta_{rt}, a_r^t)_{r,t}$ and the corresponding $\widetilde{G}$, define $\chi$ as in \S \ref{secIR}.
The arguments in \ref{lem4.5}, \ref{lem4.7}, \ref{lem4.8} still work.
The theorem follows.
\end{Proof}

\begin{Theorem}\label{thm-R1-gen}
  Let $V = \sum_{s} \mathbb{C} \alpha_s$ be an R-representation, then $V$ is a quotient of a unique IR-representation $V^\prime = \bigoplus_s \mathbb{C} \alpha_s^\prime$, say, $\pi : V^\prime \twoheadrightarrow V$. Moreover, $W$ acts on $\ker \pi$ trivially.
\end{Theorem}

\begin{Proof}
Note that we are allowing $m_{rt}$ to be $\infty$, so $\alpha_r$ and $\alpha_t$ may be proportional.
If this happens, choose $\delta_{rt} = \varrho_4$ (i.e. the geometric representation), and choose $a_r^t, a_t^r \in \mathbb{C}^\times$ such that $a_r^t \alpha_r + a_t^r \alpha_t = 0$.
Except this case, the proof is nearly the same as that of \ref{thm-R1}\ref{5.6.1}.
\end{Proof}

\smallskip\noindent\textbf{The matrix {\bfseries \itshape A}}
Let $V$ be an IR-representation.
In view of \ref{thm-R1-gen}, we wish to find sub-representations in $V$ with trivial $W$-action.
Write $S = \{s_1, \dots, s_n\}$, and use notations $P_{ij}, \delta_{ij}, \varrho_i^j$, etc in the same way as before.
Let $v = \sum_{j} x_j \alpha_j \in V$, then
\begin{align}\refstepcounter{Theorem}\label{eq-7.3}
    s_i \cdot v & = -x_i \alpha_i + \sum_{s_j \in S_1} x_j \Bigl(\alpha_j + 2 \frac{a_i^j}{a_j^i} \cos \frac{k_{ij} \uppi}{m_{ij}} \alpha_i\Bigr) + \sum_{s_j \in S_2} x_j \alpha_j \notag \\
    & \mathrel{\phantom{=}} + \sum_{s_j \in S_3} x_j \Bigl(\alpha_j + \frac{a_i^j}{a_j^i} \alpha_i \Bigr) + \sum_{s_j \in S_4} x_j \Bigl( \alpha_j + u_{ij} \frac{a_i^j}{a_j^i} \alpha_i\Bigr) \\
    & = \Big(- x_i + \sum_{s_j \in S_1} 2x_j \frac{a_i^j}{a_j^i} \cos \frac{k_{ij} \uppi}{m_{ij}} + \sum_{s_j \in S_3} x_j \frac{a_i^j}{a_j^i} + \sum_{s_j \in S_4} x_j u_{ij} \frac{a_i^j}{a_j^i} \Big) \alpha_i + \sum_{j \neq i} x_j \alpha_j, \notag
\end{align}
where
\begin{align*}
    S_1 & := \{s_j \in S \mid m_{ij} < \infty\}, \quad 1 \leq k_{ij} \leq \frac{m_{ij}}{2} \text{ such that } \delta_{ij} = \rho_{k_{ij}}, \\
    S_2 & := \{s_j \in S \mid m_{ij} = \infty, \delta_{ij} = \varrho_j^i \}, \\
    S_3 & := \{s_j \in S \mid m_{ij} = \infty, \delta_{ij} = \varrho_i^j \}, \\
    S_4 & := \{s_j \in S \mid m_{ij} = \infty, \delta_{ij} = \varrho_{z_{ij}}, z_{ij} \in \mathbb{C}\}, \\
    \quad u_{ij} & : = u(z_{ij}) \text{ chosen in \S \ref{Dm-and-graph}}.
\end{align*}
We have $S = S_1 \sqcup S_2 \sqcup S_3 \sqcup S_4 \sqcup \{s_i\}$.
Similar to \eqref{matrixA}, let $A = (A_{ij}) \in \mathrm{M}_n(\mathbb{C})$ be a square matrix with element $A_{ij}$ at $(i,j)$-position defined by (here $S_1 = S_1(i), \dots$ depend on $i$)
\begin{equation*}
    A_{ij} = \begin{cases}
      1, & \mbox{if } i = j, \\
      - \frac{a_i^j}{a_j^i} \cos \frac{k_{ij} \uppi}{m_{ij}}, & \mbox{if } s_j \in S_1, \\
      0, & \mbox{if } s_j \in S_2, \\
      - \frac{a_i^j}{2 a_j^i}, & \mbox{if } s_j \in S_3, \\
      - \frac{u_{ij} a_i^j}{2 a_j^i}, & \mbox{if } s_j \in S_4,.
    \end{cases}
\end{equation*}
From \eqref{eq-7.3}, the null space of $A$ is the subspace of $V$ with trivial $W$-action, as we saw in \S \ref{secR1}.

\smallskip\noindent\textbf{Other reducibility}
Compared with \ref{prop5.3} and \ref{prop6.3}, there are some different phenomena when considering reducibility of IR-representations in this general setting.
For example, $\varrho_r^t$ of $D_\infty$ have sub-representations with nontrivial group action.

Suppose $V$ is an IR-representation.
For convenience, let's append more information to the graph $\widetilde{G}$.
For any pair $r \neq t \in S$, if $\delta_{rt} = \varrho_r^t$, assign a direction $t \to r$ to the edge $\{r,t\}$.
Symmetrically, $r \to t$ if $\delta_{rt} = \varrho_t^r$.
For other edges $\{r,t\}$ in $\widetilde{G}$ ($3 \leq \widetilde{m}_{rt} \leq \infty$), both directions are assigned, i.e. $r \to t$ and $t \to r$.
These directions have the following meaning: if $r \to t$ then $\alpha_t$ belongs to the sub-representation generated by $\alpha_r$.

The directions generate a pre-order on $S$: say $r \dleq t$ if there is a sequence $r = s_1, s_2, \dots, s_n = t$ in $S$ such that $s_{i+1} \to s_i$ for all $i$.
The corresponding equivalence relation is denoted by $\dequiv$: say $r \dequiv t$ if $r \dleq t$ and $t \dleq r$.
Immediately we have the following fact.

\begin{Corollary}
Given an IR-representation $V$.
\begin{enumerate}
    \item If all elements in $S$ are equivalent with respect to $\dequiv$, then the action of $W$ on any sub-representation of $V$ is trivial.
    \item Let $I \subsetneq S$ be an equivalent class with respect to $\dequiv$, then $\bigoplus_{s \dleq I} \mathbb{C} \alpha_s$ is a sub-representation of $V$.
\end{enumerate}
\end{Corollary}

\smallskip\noindent\textbf{Bilinear forms}
Different from \ref{6.1}, it may happen that there is no nonzero $W$-invariant bilinear form on an IR-representation $V$ even when $\chi (\operatorname{H}_1 (\widetilde{G})) \subseteq \{\pm 1\}$.
For example, $S = \{s, t, r\}$, $m_{st} = \infty$, $m_{rt} = m_{sr} = 3$, let $\delta_{st} = \varrho_s^t$.
Suppose $\operatorname{B}$ is a $W$-invariant bilinear form on $V$.
From \ref{lem2.11}\ref{2.4.4}, it holds $\operatorname{B} (\alpha_s, \alpha_s) = 0$.
By \ref{lem-2.7}, it yields $\operatorname{B} (\alpha_r, \alpha_r) = 0$ and then $\operatorname{B} (\alpha_t, \alpha_t) = 0$.
$\operatorname{B}$ has to be zero.

If for any pair $r \neq t \in S$ it holds $\delta_{rt} \neq \varrho_r^t$ and $\varrho_t^r$, then \ref{6.1} still works.

\section{Examples: from \texorpdfstring{$\widetilde{\mathsf{A}}_n$}{affine An} to simply laced type} \label{sec-eg}
\begin{Example}
The Coxeter system $(W,S)$ of type $\widetilde{\mathsf{A}}_n$ is defined by $\langle s_0, s_1, \dots, s_n \mid s_i^2 = (s_i s_{i+1})^3 = 1, \forall i \rangle$ (regard $n+1$ as $0$).
Its Coxeter graph $G$ is
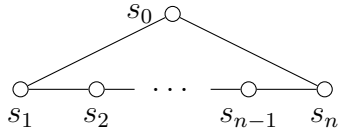
\begin{equation*}
  \begin{tikzpicture}
    \node [circle, draw, inner sep=2pt, label=below:$s_1$] (s1) at (-2,0) {};
    \node [circle, draw, inner sep=2pt, label=below:$s_2$] (s2) at (-1,0) {};
    \node [circle, draw, inner sep=2pt, label=below:$s_{n-1}$] (sn-1) at (1,0) {};
    \node [circle, draw, inner sep=2pt, label=below:$s_n$] (sn) at (2,0) {};
    \node [circle, draw, inner sep=2pt, label=left:$s_0$] (s0) at (0,1) {};
    \graph{(sn-1) -- (sn) -- (s0) -- (s1) -- (s2);};
    \draw (s2) -- (-0.5,0);
    \draw (sn-1) -- (0.5,0);
    \node (d) at (0,0) {$\dots$};
  \end{tikzpicture}
\end{equation*}

In the data defining IR-representations, for any distinct pair $i,j$, the condition $1 \leq k_{ij} \leq m_{ij}/2$ forces $k_{ij} = 1$, thus the graph $\widetilde{G} = G$.
The homology group $\operatorname{H}_1 (G)$ is isomorphic to $\mathbb{Z}$ generated by the circuit $c = (s_0, s_1, \dots, s_n, s_0)$.
Giving a character $\chi$ of $\operatorname{H}_1 (G)$ is equivalent to giving a number $x \in \mathbb{C}^\times$ and assigning $\chi(c) = x$.
Thus, all IR-representations of $W$ is parameterized by $\mathbb{C}^\times$.

Let $a_0^n = x$, and $a_n^0 = a_0^1 = a_1^0 \dots = 1$, the square matrix $A$ (of size $n + 1$, defined in \eqref{matrixA}) is
\begin{equation*}
    A = \begin{pmatrix}
        1 & - \frac{1}{2} & 0 & \dots & 0 & - \frac{x}{2} \\
        - \frac{1}{2} & 1 & - \frac{1}{2} &  & & 0 \\
        0 & - \frac{1}{2} & \ddots & \ddots &  & \vdots \\
        \vdots &  & \ddots & \ddots & - \frac{1}{2} & 0 \\
        0 & &  & - \frac{1}{2} & 1 & - \frac{1}{2} \\
        - \frac{1}{2x} & 0 & \dots & 0 & - \frac{1}{2} & 1
      \end{pmatrix}
\end{equation*}
One may use induction on $n$ to compute the determinant of $A$,
\begin{equation*}
    \det A = \frac{2 - x - x^{-1}}{2^{n+1}}.
\end{equation*}
Denote by $V_x$ the IR-representation given by $x$.
Then, by \ref{prop6.3}, $V_x$ is irreducible if and only if $x \neq 1$.
When $x = 1$, $V_x \simeq \Vg$.
\end{Example}

\begin{Proposition}\label{prop-affineAn}
Let $(W,S)$ be of type $\widetilde{\mathsf{A}}_n$.
Suppose $V$ is a nontrivial irreducible A1-representation of $W$, then either $V$ is the simple quotient of $\Vg$, or $V = V_x$ for some $x \in \mathbb{C} \setminus \{0,1\}$.
{\upshape (See also \cite[Ch.12]{Xi94} for the version of affine Hecke algebras.)}
\end{Proposition}

\begin{Proof}
First of all, $V$ is finite dimensional (see, for example, \cite{Kato83} or \cite{Xi07}).
Restricted to the dihedral group $\langle s_0, s_1 \rangle \simeq D_3$, $V$ decomposes into a direct sum of copies of $\mathds{1}, \varepsilon$ and $\rho_1$.
Let $V_-^i$ denote the subspace of $-1$-eigenvectors of $s_i$.
The condition \eqref{A1} forces that $\varepsilon$ doesn't occur in $V$, thus $V_-^0 \cap V_-^1 = 0$ and $\dim V_-^0 = \dim V_-^1$.
Similar result holds for any pair $i, i+1$.

If $V$ is nontrivial, then $V_-^0 \ne 0$.
Take a basis of $V_-^0$, say $\{\alpha_{0,j}\}_{j = 1, \dots, r}$.
Let
\begin{equation*}
  \alpha_{1,j} := s_1 \cdot \alpha_{0,j} - \alpha_{0,j}, \quad j = 1, \dots, r,
\end{equation*}
then $\{\alpha_{1,j}\}_j$ is a basis of $V_-^1$, and for any $j$, $\{\alpha_{0,j}, \alpha_{1,j}\}$ spans a representation of $\langle s_0, s_1 \rangle$ isomorphic to $\rho_1$.
Apply similar considerations inductively on $\langle s_1, s_2 \rangle, \dots, \langle s_{n-1}, s_n \rangle$, we obtain a basis $\{\alpha_{i,j}\}_j$ of $V_-^i$ ($1 \leq i \leq n$), such that $\forall 1 \leq i \leq n, 1 \leq j \leq r$, $\{\alpha_{i-1,j}, \alpha_{i,j}\}$ spans a $\rho_1$ of $\langle s_{i-1}, s_i \rangle$.

Now consider $\langle s_n, s_0 \rangle$.
Write
\begin{equation*}
  \alpha_{0,j}^\prime := s_0 \cdot \alpha_{n,j} - \alpha_{n,j}, \quad j = 1, \dots, r,
\end{equation*}
then $\{\alpha_{n,j}, \alpha_{0,j}^\prime\}$ spans a $\rho_1$ of $\langle s_n, s_0 \rangle$ for any $j$, and $\{\alpha_{0,j}^\prime\}_j$ is a basis of $V_-^0$ as well.
Thus there is a matrix $X \in \operatorname{GL}_r(\mathbb{C})$ such that
\begin{equation*}
    (\alpha_{0,1}^\prime, \dots, \alpha_{0,r}^\prime) = (\alpha_{0,1}, \dots, \alpha_{0,r}) \cdot X.
\end{equation*}
Let $v = (x_1, \dots, x_r)^\mathrm{T} \in \mathbb{C}^r$ be an eigenvector of $X$ of with eigenvalue $x \in \mathbb{C}^\times$.
Write
\begin{equation*}
    \alpha_i := \sum_{1\leq j \leq r} x_j \alpha_{i,j}, \quad 0 \leq i \leq n, \quad \text{and } \alpha_0^\prime := \sum_{1\leq j \leq r} x_j \alpha_{i,j}^\prime,
\end{equation*}
then for $1 \leq i \leq n$, $\{\alpha_{i-1}, \alpha_i\}$ span a $\rho_1$ for $\langle s_{i-1}, s_i \rangle$.
At last, $\{\alpha_n, \alpha_0^\prime\}$ span a $\rho_1$ for $\langle s_n, s_0 \rangle$.
We have
\begin{equation*}
    \alpha_0^\prime = (\alpha_{0,1}^\prime, \dots, \alpha_{0,r}^\prime) \cdot v = (\alpha_{0,1}, \dots, \alpha_{0,r}) \cdot X \cdot v = x (\alpha_{0,1}, \dots, \alpha_{0,r}) \cdot v = x \alpha_0.
\end{equation*}

Consider the subspace $V^\prime$ spanned by $\{\alpha_0, \alpha_1, \dots, \alpha_n\}$, we claim that $V^\prime$ forms a sub-representation satisfying \eqref{R}.
For any pair $i, i^\prime$ such that $s_i, s_{i^\prime}$ are not adjacent in $G$, we have $s_i s_{i^\prime} = s_{i^\prime} s_i$, and hence $s_i \cdot \alpha_{i^\prime} \in V_-^{i^\prime}$.
Then $s_i \cdot \alpha_{i^\prime} - \alpha_{i^\prime} \in V_-^i \cap V_-^{i^\prime}$, thus $s_i \cdot \alpha_{i^\prime} = \alpha_{i^\prime}$.
Now we can see that for any $i$, $s_i$ leaves $V^\prime$ stable, and has only one $-1$-eigenvector in $V^\prime$ (up to scalars).
Thus, $V^\prime$ is an R-representation.
By the construction of $V^\prime$ and the classification of R-representations, $V^\prime$ is isomorphic to $V_x$ if $x \neq 1$, to $\Vg$ or its simple quotient if $x = 1$.
The proposition is proved.
\end{Proof}

This proposition says, for type $\widetilde{\mathsf{A}}_n$ any nontrivial irreducible A1-representation is an R-representations.
The proof can be generalized a little bit, as follows.

\begin{Theorem} \label{8.3}
Let $(W,S)$ be an irreducible Coxeter system with simply laced Coxeter graph $G$.
Suppose the rank of $\operatorname{H}_1(G)$ is $0$ or $1$, then any nontrivial irreducible A1-representation is an R-representation.
In particular, they are finite dimensional.
If $\operatorname{H}_1(\Gamma)$ has rank $0$, such representation is unique.
\end{Theorem}

\begin{Proof}
Let $V$ be an irreducible A1-representation of $W$, denote by $V_-^s$ the subspace of $-1$-eigenvectors of $s$.
If $\operatorname{H}_1(G) \simeq \mathbb{Z}$, let $c$ be a circuit generating $\operatorname{H}_1(G)$, and $s$ be a vertex in $c$.
Similar to the proof of \ref{prop-affineAn}, the circuit $c$ gives a linear automorphism $X$ of $V_-^s$.
Now $V_-^s$ becomes a module of the Laurent polynomial ring $\mathbb{C}[X^{\pm 1}]$.
By Hilbert's Nullstellensatz, any simple module of $\mathbb{C}[X^{\pm 1}]$ is one dimensional.
Therefore, if $\dim V_-^s > 1$, then there is a proper submodule $V_0 \subsetneq V_-^s$.
This subspace $V_0$ generates a proper sub-representation of $W$ in $V$, contradicting the irreducibility of $V$.
Hence $\dim V_-^s = 1$, and $V_-^s$ generates an R-representation.

If $G$ is a tree, fixing $s \in S$, any nonzero vector in $V_-^s$ generates a sub-representation, which is an R-representation.
Such a representation is unique according to \ref{thm-R1}.
\end{Proof}

Motivated by this proof, for a Coxeter graph $G$ with at least two circuits, we may use an infinite dimensional irreducible representation of the fundamental group of $G$ (which is a non-abelian free group) to construct an infinite dimensional irreducible representation of the corresponding Coxeter group.
We will discuss this in a subsequent paper, see \cite{Hu22}.

\begin{Example}\label{7.2}
This example is a supplement to \ref{rmk5.7}\ref{5.7.2}.
When $\widetilde{G}$ is not connected, the quotients of an IR-representation by two different subspace may be isomorphic.
Suppose $\widetilde{G}$ is of type $\widetilde{\mathsf{A}}_2 \times \widetilde{\mathsf{A}}_2$, i.e.
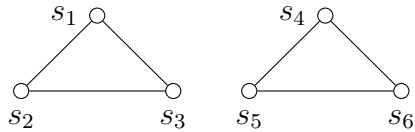
\begin{equation*}
    \begin{tikzpicture}
      \node [circle, draw, inner sep=2pt, label=left:$s_1$] (s1) at (-1.5,1) {};
      \node [circle, draw, inner sep=2pt, label=below:$s_2$] (s2) at (-2.5,0) {};
      \node [circle, draw, inner sep=2pt, label=below:$s_3$] (s3) at (-0.5,0) {};
      \graph{(s3) -- (s1) -- (s2) -- (s3);};
      \node [circle, draw, inner sep=2pt, label=left:$s_4$] (s4) at (1.5,1) {};
      \node [circle, draw, inner sep=2pt, label=below:$s_5$] (s5) at (0.5,0) {};
      \node [circle, draw, inner sep=2pt, label=below:$s_6$] (s6) at (2.5,0) {};
      \graph{(s4) -- (s5) -- (s6) -- (s4);};
    \end{tikzpicture}
\end{equation*}
and let $V$ be the corresponding geometric representation.
Then $V$ is isomorphic to the direct sum of two copies of the geometric representations of $\widetilde{\mathsf{A}}_2$.
Let $v_1 = \alpha_1 + \alpha_2 + \alpha_3$, $v_2 = \alpha_4 + \alpha_5 + \alpha_6$, then $V_0 := \langle v_1, v_2 \rangle$ is the maximal sub-representation in $V$ with trivial group action.
The two quotient $V/\langle v_1 \rangle$ and $V/\langle v_2 \rangle$ are clearly non-isomorphic.
However, let $V_1 = V/\langle v_1 + v_2 \rangle$, $V_2 = V/\langle v_1 + 2v_2 \rangle$, then $\alpha_i \mapsto \alpha_i, i = 1,2,3$, $\alpha_j \mapsto 2 \alpha_j, j = 4,5,6$ defines an isomorphism from $V_1$ to $V_2$.
\end{Example}

\pdfbookmark[1]{References}{references}
\bibliography{reps-and-h1}

\bigskip \normalsize
\textsc{Academy of Mathematics and Systems Science, Chinese Academy of Sciences, Beijing, China}

\textit{Email address:} \texttt{huhongsheng16@mails.ucas.ac.cn}
\end{document}